\documentclass[preprint, 11pt]{elsarticle}

\usepackage[T1]{fontenc}
\usepackage[utf8]{inputenc}
\usepackage[hmargin=0.9in]{geometry}
\usepackage[babel,german=quotes]{csquotes}
\usepackage{subfig}
\usepackage{graphicx}
\usepackage[colorlinks=true,citecolor=blue,urlcolor=blue]{hyperref}
\usepackage{caption}
\usepackage{amsmath}          
\usepackage{amsthm}           
\usepackage{amssymb}          
\usepackage{bm, amsfonts}     
\usepackage{xcolor}
\usepackage{fixmath}
\usepackage{tikz}
\usepackage{bm}
\usepackage{makecell}

\graphicspath{{./pictures/small/}}
\newdefinition{cor}{Corollary}
\newdefinition{example}{Example}
\newproof{pf}{Proof}

\newcommand {\dx} {\,{\rm d}{\mathbf x}}

\newcommand {\ds} {\,{\rm d}{\mathrm s}}

  \newcommand{\R}{\mathbb{R}}
  \newdefinition{rmk}{Remark}
  
  \newcommand{\pd}[2]{\frac{\partial #1}{\partial #2}}

  \newcommand{\td}[2]{\frac{\mathrm d #1}{\mathrm d #2}}

\newcommand{\beq}{\begin{equation}}
\newcommand{\eeq}{\end{equation}}

\newcommand{\Real}{\mathbb{R}}

\makeatletter
\def\ps@pprintTitle{%
  \let\@oddhead\@empty
  \let\@evenhead\@empty
  \def\@oddfoot{
    \footnotesize\itshape
    \hfill\today
  }%
  \let\@evenfoot\@oddfoot}
\makeatother

\begin{document}

\begin{frontmatter} 
  \title{Optimal control using flux potentials: A way to construct bound-preserving finite element schemes for conservation laws}

\author{Falko Ruppenthal}
\ead{falko.ruppenthal@math.tu-dortmund.de}

\author{Dmitri Kuzmin\corref{cor1}}
\ead{kuzmin@math.uni-dortmund.de}

\cortext[cor1]{Corresponding author}

\address{Institute of Applied Mathematics (LS III), TU Dortmund University\\ Vogelpothsweg 87,
  D-44227 Dortmund, Germany}

\journal{Journal of Computational Physics}

\begin{abstract}  
  To ensure preservation of local or global bounds for numerical solutions
  of conservation laws, we constrain a baseline finite element discretization
  using optimization-based (OB) flux correction. The main novelty of
  the proposed
  methodology lies in the use of flux potentials as control variables and
  targets of inequality-constrained optimization problems for numerical
  fluxes.  In contrast to optimal control via
  general source terms, the discrete conservation property of
  flux-corrected finite element approximations
  is guaranteed without the need to impose
  additional equality constraints.
Since the number of flux potentials is less than the number
  of fluxes in the multidimensional case, the potential-based version
  of optimal flux control involves fewer unknowns than direct
  calculation of optimal fluxes.  We show that the feasible set
  of a potential-state potential-target (PP) optimization problem
  is nonempty and choose a primal-dual Newton method for
  calculating the optimal flux potentials. The results of
  numerical studies for linear advection and anisotropic 
  diffusion problems in 2D demonstrate the superiority of the new OB-PP
  algorithms to closed-form flux limiting under worst-case assumptions.

\end{abstract}
\end{frontmatter}

\section{Introduction}

In many applications of practical interest, it is essential to
guarantee that numerical approximations to a scalar conserved quantity
$u$ attain values in a range $[u^{\min},u^{\max}]$ of physically admissible
states. For example, concentrations and volume
fractions
must stay between $u^{\min}=0$ and $u^{\max}=1$.
Because of modeling errors, even exact solutions of the governing
equations can sometimes violate such constraints.
A physics-compatible numerical method should ensure preservation of
global bounds and the discrete conservation property, while
keeping consistency errors as small as possible. Moreover,
numerical solutions may be required to satisfy
local maximum principles and stay free of spurious oscillations.

A very general framework for enforcing inequality constraints
in numerical schemes for conservation laws is based on the
concept of algebraic flux correction (AFC) \cite{afc_analysis,afc1,convex}.
Given a high-order baseline discretization,
an AFC scheme modifies it using numerical fluxes associated with
a graph Laplacian operator. The validity of relevant maximum
principles is enforced using flux limiters. A typical limiting
technique is derived by formulating an inequality-constrained
optimization problem and making worst-case assumptions to
obtain a closed-form formula for limited fluxes or correction
factors. Many limiters of this kind are based on ideas introduced
in the work of Boris and Book \cite{fct1}, Van Leer \cite{vl1,vl2}, 
Zalesak \cite{zalesak79}, and Harten \cite{harten1,harten2}.
If necessary, the accuracy of flux-corrected approximations
can be improved by using
iterative or optimization-based limiting procedures. For example,
iterative methods for calculating 
physics-aware flux approximations were proposed in
\cite{CH,fctmax}. 

The development of optimization-based (OB) finite element methods 
for conservation laws was greatly advanced by the recent work of Bochev
et al. \cite{Bochev2020,Bochev2014,Bochev2012} who used different
kinds of optimal control to ensure preservation of local bounds.
Their flux-state flux-target (FF) optimization method \cite{Bochev2020}
has the structure of an AFC scheme. The equivalence to a flux-corrected
transport (FCT) algorithm was shown in \cite{Bochev2020}
for a simplified quadratic
programming problem with box constraints. This interesting 
relationship means that flux limiters can be interpreted as
approximate solvers for optimization problems. In particular,
the existence of a bound-preserving
AFC approximation implies that the feasible
set is nonempty and provides a good initial guess for iterative
optimization.

If fluxes are used as optimization variables, as in the OB-FF
method, the size of the problem depends on the number
of edges in the sparsity graph of the finite element matrix. To
reduce the number of unknowns and choose a more intuitive target
for AFC, we formulate inequality-constrained optimization
problems using {\it flux potentials} in the present paper.
The new optimization variables can be interpreted as approximations
to nodal time derivatives, and the size of the resulting 
potential-state potential-target (PP) optimization problem
is proportional to the number of nodes (rather than edges).
The optimal fluxes,  defined in terms of the 
potentials, can be applied after the discretization
in space and time or at the level of spatial semi-discretization.
The latter option makes it possible to construct OB-PP counterparts
of the monolithic AFC schemes developed in \cite{convex,entropyHO,timelim}.

In the next two sections, we discretize a generic conservation
law using the continuous Galerkin method and review the basics of
algebraic flux correction. In Section \ref{sec:limiting}, we
discuss closed-form flux limiting and its connection to
optimal control of OB-FF type. The new OB-PP approach is 
introduced in Section \ref{sec:opt}, where we define the flux
potentials, formulate the inequality constraints, and show
that the feasible set is nonempty. A primal-dual Newton's
method for solving OB-PP optimization problems is outlined in
Section \ref{sec:solver}. The numerical examples of Section 
\ref{sec:num} demonstrate the potential benefits of optimal
control for numerical approximations to linear advection and
anisotropic diffusion problems. We conclude this paper with
a discussion of the main results and open problems in
Section \ref{sec:conclusions}.

\section{Baseline discretization}
\label{sec:baseline}
  
Let $u(\mathbf{x},t)$ be a scalar conserved
quantity depending on the space location $\mathbf{x}\in
\R^d,\ d\in\{1,2,3\}$ and time instant $t\ge 0$. Restricting
our attention to a bounded domain $\Omega\subset\R^d$
with  Lipschitz boundary $\Gamma=\partial\Omega$, we
consider a general initial-boundary value problem of the form
\begin{subequations}\label{ibvp}
\begin{alignat}{3}
  \pd{u}{t}+\nabla\cdot\mathbf{f}(u)&=\nabla\cdot(\mathcal D\nabla u)
  &&\qquad\mbox{in}\ \Omega\times\R_+,
 \label{ibvp-pde}\\
  u&=u_D&&\qquad\mbox{on}\ \Gamma_D\times\R_+,\label{ibvp-bc}\\
  u&=u_0&&\qquad\mbox{on}\ \Omega\times\{0\},\label{ibvp-ic}
\end{alignat}
\end{subequations}
where $\mathbf{f}(u)\in\R^d$ is an inviscid flux and $\mathcal D
\in\R^{d\times d}$
is a symmetric positive semidefinite tensor.
The Dirichlet boundary data $u_D$ is prescribed on 
$\Gamma_D=\Gamma$ in the case $\mathcal D\ne 0$. For hyperbolic problems,
$\Gamma_D=\{\mathbf{x}\in\Gamma\,:\,\mathbf{f}'(u)\cdot\mathbf{n}
(\mathbf{x})<0\}$,
where $\mathbf{n}$ is the unit outward normal.
We are also interested in steady-state solutions of
problem \eqref{ibvp}, i.e., in solutions of 
the boundary value problem 
\begin{subequations}\label{bvp}
\begin{alignat}{3}
  \nabla\cdot\mathbf{f}(u)&=\nabla\cdot(\mathcal D\nabla u)
  &&\qquad\mbox{in}\ \Omega,
 \label{bvp-pde}\\
  u&=u_D&&\qquad\mbox{on}\ \Gamma_D.\label{bvp-bc}
\end{alignat}
\end{subequations}

We discretize \eqref{ibvp-pde} in space using a numerical scheme
that yields an algebraic system of the form
\beq\label{semi-c}
M_C\td{u}{t}=r(u).
\eeq
By abuse of notation, $u=(u_1,\ldots,u_{N_h})^T$ denotes the vector
of $N_h$ discrete unknowns in formulas involving matrices rather than
differential operators. The subscript $C$ is used for consistent
mass matrices of (differential-)algebraic systems corresponding to a
baseline scheme, such as the continuous Galerkin discretization
with linear ($\mathbb{P}_1$) or multilinear ($\mathbb{Q}_1$)
finite elements. For this particular method, a weak form 
of  \eqref{ibvp-pde}, \eqref{ibvp-bc} yields
$M_C=(m_{ij})_{i,j=1}^{N_h}$ and $r(u)=(r_{i})_{i=1}^{N_h}$
with entries
$$
m_{ij}=\sum_{e=1}^{E_h}\int_{K_e}\varphi_i\varphi_j\dx,
$$
\beq\label{r-def}
\qquad r_{i}=
\sum_{e=1}^{E_h}\left[\lambda
\int_{K_e\cap\Gamma_D}
\varphi_i(u_D-u_h)\ds-
  \int_{K_e}\varphi_i\nabla\cdot\mathbf{f}(u_h)\dx-
\int_{K_e}
\nabla\varphi_i\cdot(\mathcal D\nabla u_h)\dx
  \right],
\eeq
where $\lambda>0$ is a sufficiently large penalty parameter,
$K_1,\ldots,K_{E_h}$ are elements of a 
computational mesh and $\varphi_1,\ldots,\varphi_{N_h}$ are the basis
functions of the Lagrange finite element approximation
$$
u_h(\mathbf{x},t)=\sum_{j=1}^{N_h}u_j(t)\varphi_j(\mathbf{x}).
$$
The flux function $\mathbf{f}(u_h)$ is sometimes approximated by the
{\it group finite element} interpolant \cite{group,afc1}
$$
\mathbf{f}_h(\mathbf{x},t)=\sum_{j=1}^{N_h}\mathbf{f}(u_j(t))
\varphi_j(\mathbf{x})
$$
and its contribution to $r(u)$ is stabilized using additional
terms. In the numerical experiments
of Section~\ref{sec:num}, we use Taylor-Galerkin
stabilization \cite{quartapelle,selmin_fef93}
for linear advection problems.

To discuss the modifications that are needed to satisfy
maximum principles for general discretizations
and conservation laws, we abstain from
giving further details of our baseline method so far.

\section{Algebraic flux correction}
\label{sec:afc}

To enforce the validity of relevant inequality constraints using
the algebraic flux correction (AFC) methodology \cite{afc1}, we
replace \eqref{semi-c} by the system of
ordinary differential equations
\beq\label{semi-l}
M_L\td{u}{t}=r(u)+g(u),
\eeq
where $M_L=\mathrm{diag}(\delta_{ij}m_i)$ is the
lumped mass matrix with positive diagonal entries
$$m_i=\sum_{j=1}^{N_h}m_{ij}>0.$$
The vector $g(u)=(g_i)_{i=1}^{N_h}$ defines a particular semi-discrete
scheme and can be interpreted as source term control for
inequality-constrained optimization. Note that systems  \eqref{semi-c} and
  \eqref{semi-l} are equivalent for
  $g(u)=(M_L-M_C)\dot u$, where $\dot u$ is the solution of the
linear system $M_C\dot u=r(u)$.

To preserve the discrete conservation property, the control 
term $g(u)$ must satisfy $\sum_{i=1}^{N_h}g_i=0.$ 
This zero sum condition implies
existence of numerical fluxes $g_{ij}=-g_{ji}$ such that
$$
g_i=\sum_{j\in\mathcal N_i\backslash\{i\}}g_{ij},
$$
where $\mathcal N_i$ is the  stencil of node $i$. For our finite element
scheme
$
\mathcal N_i=\{j\in\{1,\ldots N_h\}\,:\,m_{ij}> 0\}
$ is the integer set containing the index $i$ and the indices of
nearest neighbors of node $i$.

The fluxes corresponding to the high-order
target $g^H(u)=(M_L-M_C)\dot u$ are given by \cite{afc1,convex}
\beq\label{g-target}
g_{ij}^H=m_{ij}(\dot u_i-\dot u_j),\qquad \dot u=M_C^{-1}r(u).
\eeq
Indeed, the $i$th component of  $g^H(u)$ admits the following
decomposition:
$$
g_i^H=m_i\dot u_i-\sum_{j\in\mathcal N_i}m_{ij}\dot u_j
=\sum_{j\in\mathcal N_i}m_{ij}\dot u_i
-\sum_{j\in\mathcal N_i}m_{ij}\dot u_j=
\sum_{j\in\mathcal N_i\backslash\{i\}}m_{ij}(\dot u_i-\dot u_j).
$$

Let us discretize the AFC system \eqref{semi-l} in time using an explicit or
implicit method that yields
\beq\label{full-c}
M_Lu^{n+1}=M_Lu^n+\Delta t(r+g),
\eeq
where $\Delta t>0$ is a constant
time step and $u^n\approx
u(t^n)$ for
$t^n=n\Delta t$. The number of intermediate stages (if any)
and the way to compute
$r+g$ depend on the particular space-time discretization.

A numerical scheme  is called {\em bound preserving} (BP)
if discrete maximum principles of the form
\beq\label{lmp}
u^{\min}\le u_i^{\min}(t)\le u_i(t) \le u_i^{\max}(t)\le u^{\max}
\eeq
hold for each node and time level.
As shown in \cite{timelim}, the BP property of \eqref{full-c}
is guaranteed
for any strong stability preserving (SSP) Runge-Kutta time
integrator if \eqref{semi-l} can be written as
\beq\label{pp-mono}
m_i\td{u_i}{t}=c_i(u_i^*-u_i),\qquad i=1,\ldots,N_h,
\eeq
where $c_i>0$ is bounded and $u_i^*$ is an intermediate state such that
\beq\label{pp-mono2}
u_j\in [u_i^{\min},u_i^{\max}]\quad\forall
j\in\mathcal N_i\quad\Rightarrow\quad
u_i^*\in [u_i^{\min},u_i^{\max}].
\eeq
\begin{rmk}
For an explicit SSP-RK method with forward Euler stages of the form
$$
u_i^{\rm SSP}=u_i+\frac{\Delta t}{m_i}c_i(u_i^*-u_i)
$$
and a time step satisfying 
$\Delta t c_i\le m_i$, the result $u_i^{\rm SSP}=(1-\Delta t c_i/m_i)u_i
+\Delta t c_i/m_i)u_i^*$ is a convex combination of the states
$u_i$ and $u_i^*$.
It follows that $u_i^{\rm SSP}\in [u_i^{\min},u_i^{\max}]$.
\end{rmk}
The objective of algebraic flux correction is to enforce
the inequality constraints for $u_i^*$ or $u_i^n$
using sufficiently accurate approximations $g_{ij}^*$
to the fluxes $g_{ij}$ of the target discretization.

\section{Flux limiting and optimization}
\label{sec:limiting}

Many generalizations of classical flux-corrected transport
(FCT) algorithms \cite{fct1,zalesak79} and total variation
diminishing (TVD) methods \cite{harten1,harten2}, combine
\eqref{semi-c} with a low-order BP semi-discretization
$$
M_L\td{u}{t}=r(u)+Du,
$$
where $D=(d_{ij})_{i,j=1}^{N_h}$ is an artificial
diffusion
operator such that $d_{ij}\ge 0$ for $j\ne i$. If $D$ is
a graph Laplacian (i.e., a symmetric
matrix with zero row and column sums), the
corresponding low-order
AFC control $g^L(u)=Du$ admits a conservative
decomposition into the diffusive fluxes \cite{afc1}
\beq\label{g-low}
g_{ij}^L=d_{ij}(u_j-u_i).
\eeq
Such fluxes are used to enforce the BP
property in the vicinity of discontinuities and steep gradients.
In smooth regions, the Galerkin approximation to a hyperbolic
or advection-dominated conservation law may be stabilized using
a high-order dissipative flux $g_{ij}^S$. The use of high-order
stabilization in numerical fluxes of AFC approximations makes it
possible to avoid spurious ripples, achieve optimal convergence
behavior, and ensure entropy stability for nonlinear problems
\cite{Guermond2014,convex,entropyHO}.

In nonlinear high-resolution schemes, a convex combination of
$g_{ij}^L$ and $g_{ij}=g_{ij}^H+g_{ij}^S$ is defined by
\beq
g_{ij}^*=g_{ij}^L+\alpha_{ij}f_{ij}
=\alpha_{ij}g_{ij}+(1-\alpha_{ij})d_{ij}(u_j-u_i),
\eeq
where $\alpha_{ij}\in[0,1]$ is a correction factor
satisfying the symmetry condition $\alpha_{ji}=\alpha_{ij}$ and
\beq\label{f-target}
f_{ij}=m_{ij}(\dot u_i-\dot u_j)+d_{ij}(u_i-u_j)+g_{ij}^S
\eeq
is a raw antidiffusive flux such that $f_{ji}=-f_{ij}$.
An algorithm for calculating $\alpha_{ij}$ is
called a flux limiter.

Since the solution $u^{n+1}$ of the fully
discrete problem \eqref{full-c} and the
auxiliary states $u_i^*$ of the space discretization \eqref{pp-mono}
depend on the sum of limited
antidiffusive fluxes $f_{ij}^*=\alpha_{ij}f_{ij}$, the BP property 
\eqref{lmp} of an AFC scheme
can typically be shown under flux constraints of the form
\cite{afc1,convex,zalesak79}
$$
f_i^{\min}\le \sum_{j\in\mathcal N_i\backslash\{i\}}\alpha_{ij}f_{ij}
\le f_i^{\max}
$$
and CFL-like time step restrictions. To avoid solution of a global
inequality-constrained optimization problem and derive a closed-form
expression for $\alpha_{ij}$, flux limiting is usually performed
under worst-case assumptions. Examples of AFC schemes equipped with
closed-form flux limiters include a family of multidimensional
FCT algorithms \cite{Guermond2014,afc1,zalesak79} and
the monolithic convex limiting (MCL) strategy proposed in 
\cite{convex}. We use FCT and MCL in some numerical
experiments of Section \ref{sec:num}.

As shown by Bochev et al. \cite{Bochev2020,Bochev2014,Bochev2012}
in the context of remap and transport algorithms, closed-form flux
limiters approximate the exact solution of an inequality-constrained
optimization problem for $\alpha_{ij}$ by the exact solution of a
simplified optimization problem with box constraints. The lack
of optimality can be cured by using iterative flux correction
\cite{CH,afc1,fctmax} or optimization-based (OB) limiting
\cite{Bochev2020,Bochev2014}. The flux-state flux-target
(FF) algorithm presented in \cite{Bochev2020} is an
AFC scheme that determines the optimal antidiffusive fluxes
$f_{ij}^*=\alpha_{ij}f_{ij}$ by solving the quadratic programming
(QP) problem
\beq\label{ff-obj}
\mbox{minimize}\ \sum_{i=1}^{N_h}\sum_{j\in\mathcal N_i\backslash\{i\}}
(f_{ij}^*-f_{ij})^2\qquad \mbox{subject to}
\eeq
\beq\label{ff-constr}
f_{ji}^*=-f_{ij}^*,\qquad 
f_i^{\min}\le \sum_{j\in\mathcal N_i\backslash\{i\}}f_{ij}^*
\le f_i^{\max}.
\eeq
If $f_i^{\min}\le 0$ and $f_i^{\max}\ge 0$ for $i=1,\ldots,N_h$,
then the feasible set of the QP problem is nonempty  because
the fluxes $f_{ij}^*=0$ satisfy conditions \eqref{ff-constr}. By default,
the target flux $f_{ij}$ is defined by \eqref{f-target}.
An initial guess for an iterative optimization procedure
can be calculated using an AFC scheme with a closed-form
limiter of FCT or MCL type. In many cases, further iterations bring about
just marginal improvements. However,  the need for iterative limiting
arises, e.g., in situations when calculation of (nearly) optimal
fluxes / correction factors is required to minimize consistency
errors \cite{CH,fctmax}.

\section{Optimal control and flux potentials}
\label{sec:opt}

The OB-FF algorithm formulated in Section \ref{sec:limiting}
is a PDE-constrained optimization method of
discretize-then-optimize type with the state equation
\eqref{full-c} and flux variables $f_{ij}^*$. To avoid formal
dependence of the flux target $f_{ij}$ on the artificial
diffusion coefficient $d_{ij}$, the FF optimization problem
can be formulated in terms of the flux variables $g_{ij}^*$ and
target fluxes $g_{ij}=g_{ij}^H+g_{ij}^S$ as follows:
\beq
\mbox{minimize}\ \sum_{i=1}^{N_h}\sum_{j\in\mathcal N_i\backslash\{i\}}
(g_{ij}^*-g_{ij})^2\qquad \mbox{subject to}
\eeq
\beq
g_{ji}^*=-g_{ij}^*,\qquad u_i^{\min}\le u_i^*=
u_i^n+\frac{\Delta t}{m_i}\left(r_i+
 \sum_{j\in\mathcal N_i\backslash\{i\}}g_{ij}^*\right)
\le u_i^{\max}.
\eeq
To reduce the number of control variables, we introduce a new
kind of optimal flux control  in this section. Mimicking
the definition \eqref{g-target} of 
$g_{ij}^H=m_{ij}(\dot u_i-\dot u_j)$,
we express the flux variables $g_{ij}^*$ in terms of optimal 
{\it flux potentials} $\dot u_i^*$ that can be interpreted
as modified time derivatives. Using this representation,
we seek the best approximation $\dot u_i^*$ to a given potential
target $\dot u_i$ such that the inequality constraints~\eqref{lmp}
hold for \eqref{full-c} with $g$ assembled from
$g_{ij}^*=m_{ij}(\dot u_i^*-\dot u_j^*)$. The resulting OB
method can be classified as a potential-state potential-target
(PP) optimization algorithm.

Let the baseline space-time
discretization of \eqref{ibvp-pde} be given by an algebraic system of the form 
\beq\label{ofc-sys}
M_Cu^T=M_Cu^n+\Delta t r\quad\Leftrightarrow\quad
M_Lu^T=M_Lu^n+\Delta t r-(M_C-M_L)(u^T-u^n),
\eeq
where $r$ may include optional high-order stabilization. Choosing
the target potential
$$
\dot u^T=\frac{u^T-u^n}{\Delta t}=M_C^{-1}r,
$$
we calculate the optimal flux potentials $\dot u_i^*$ by solving
the PP optimization problem 
\beq\label{pp-obj}
\mbox{minimize}\ \frac{1}{2} (\dot u^* - \dot u^T)^\top M_C (\dot u^* - \dot u^T)
 + \frac{\mu}{2} (\dot u^*)^\top (M_L-M_C) \dot u^*
 \qquad \mbox{subject to}
\eeq
\beq\label{pp-constr}
\frac{m_i}{\Delta t}(u_i^{\min}-\tilde u_i^T)
\le \sum_{j\in\mathcal N_i\backslash\{i\}}m_{ij}(\dot u_i^*-\dot u_j^*)
\le \frac{m_i}{\Delta t}(u_i^{\max}-\tilde u_i^T),
\eeq
where $\tilde u_i^T=u_i^n+\Delta tr_i/m_i$ is a lumped-mass approximation to $u_i^T$ defined by the solution of \eqref{ofc-sys}.
The first term in the
definition \eqref{pp-obj}
of the objective function is
$\frac12\|\dot u_h^* - \dot u_h^T\|_{L^2(\Omega)}^2$. The
second one is introduced for
stabilization purposes and
can be controlled using the parameter $\mu\ge 0$.
The so-defined optimization problem has the same
structure as optimal control approaches with elliptic
operators and pointwise constraints \cite{roesch}.
Indeed, the matrix of our algebraic
stabilization term $g=(M_L-M_C)\dot u^*$
has the properties of
a discrete Laplacian operator that fits into the
framework developed in \cite{roesch}.

\begin{rmk}
Note that conditions \eqref{pp-constr} imply the validity of the discrete maximum
principle
\beq\label{pp-constr2}
u_i^{\min}\le u_i^*=\tilde u_i^T+\frac{\Delta t}{m_i}
\sum_{j\in\mathcal N_i\backslash\{i\}}m_{ij}(\dot u_i^*-\dot u_j^*)\le u_i^{\max}
\eeq
and  $u_i^*=u_i^T$ is obtained if the PP algorithm produces
$\dot u_j^*=\dot u_j^T$ for all $j\in\mathcal N_i$. Indeed, we
have
$$
\tilde u_i^T+\sum_{j\in\mathcal N_i\backslash\{i\}}m_{ij}(\dot u_i^T-\dot u_j^T)
=u_i^n+\frac{1}{m_i}\left(\Delta t
  r_i+\sum_{j\in\mathcal N_i\backslash\{i\}}m_{ij}[(u_i^T-u_j^T)
    -(u_i^n-u_j^n)]
  \right)=u_i^T.$$
\end{rmk}

\begin{rmk}
  The linear system corresponding to an explicit
  baseline discretization of the form \eqref{ofc-sys}
  can be solved efficiently using a few iterations of the deferred
  correction method
  $$
  M_Lu^{(m)}=M_Lu^n+\Delta tr-(M_C-M_L)(u^{(m-1)}-u^n),\qquad
  m=1,\ldots,M
  $$
  with the initial guess $u^{(0)}=\tilde u^T=
  u^n+\Delta tM_L^{-1}r$ and the final
  result $u^T:=u^{(M)}$. This approach corresponds to
  approximation of $M_C^{-1}$ by a truncated Neumann
  series; see \cite{Guermond2014,quartapelle} for details.
\end{rmk}

Let us now show that the feasible set of the PP optimization problem
\eqref{pp-obj},\eqref{pp-constr}
is nonempty, i.e., that there exists a vector $\dot u^B$ of backup
potentials such that conditions
\eqref{pp-constr} and the equivalent inequality constraints
\eqref{pp-constr2} hold for $\dot u^*=\dot u^B$. To that end,
we define $u^B$ as an auxiliary solution corresponding
to a BP space discretization of the form \eqref{pp-mono}
with the right-hand side $r_i^B(u)$ such that
\beq\label{zerosum}
\sum_{i=1}^{N_h}(r_i^B-r_i)=0.
\eeq
The validity of this zero sum condition implies that
a solution $\dot u^B$ of the linear system
$$
(M_L-M_C)\dot u^B=r^B-r
$$
with the symmetric positive semidefinite
graph Laplacian $M_L-M_C$ exists and is
unique up to a constant. The corresponding
flux variables $g_{ij}^B=m_{ij}(\dot u_i^B-\dot u_j^B)$
are defined uniquely. The use of 
$\dot u^*=\dot u^B$ in \eqref{pp-constr2} yields
$u_i^*=\tilde u_i^T+\Delta t(r_i^B-r_i)/m_i=
u_i^n+\Delta t r_i^B/m_i=u_i^B$.

It remains to construct $r^B(u)$ such that 
condition \eqref{zerosum} holds and
$u_i^B\in[u_i^{\min},u_i^{\max}]$,
perhaps under time step restrictions. Using the
residual distribution procedure proposed in
\cite{RD-BFCT}, we define
$$
r_i^B=\frac{\omega_i\rho}{\sum_{j=1}^{N_h}\omega_j},\quad
   \qquad \rho=\sum_{i=1}^{N_h}r_i,\qquad
  \omega_i=\begin{cases}
  u_i^{\max}-u_i^n & \quad \mbox{if}\  \rho>0,\\
  0 & \quad \mbox{if}\  \rho=0,\\
  u_i^{\min}-u_i^n & \quad \mbox{if}\  \rho<0.
  \end{cases}
  $$
  Note that the residual components
  $r_i^B$ satisfy \eqref{zerosum} and
  $r_i^B=c^B(u_i^*-u_i^n)$, where
  $$c^B=\frac{\rho}{\sum_{j=1}^{N_h}\omega_j}\ge 0,\qquad
  u_i^*\in\{u_i^{\min},u_i^{\max}\}.$$
  Hence, the BP property of the feasible approximation
  $u_i^B=u_i^n+\Delta t r_i^B/m_i$ is guaranteed by
\eqref{pp-mono2}, at least for time steps satisfying
the {\em a posteriori} CFL-like condition $\Delta t c^B\le m_i$; cf.
\cite{RD-BFCT}.

  For $r_i$ defined by \eqref{r-def}, we have $\rho=
  \sum_{e=1}^{E_h}\left[\lambda\int_{K_e\cap\Gamma_D}(u_D-u_h)\ds
   - \int_{K_e\cap\Gamma}\mathbf{f}(u_h)\cdot\mathbf{n}\ds
   \right]$. The addition of optional high-order stabilization
  terms does not change the value of $\rho$. Hence, the coefficient
  $c^B\ge 0$ of the ``CFL'' constraint is independent of these terms and
  of  the physical diffusion tensor  $\mathcal D$.

\begin{rmk}
  Implicit schemes with closed-form limiters may ensure the
  BP property under milder time step restrictions or unconditionally.
  If such a scheme exists for the given problem, the
  corresponding flux potential $\dot u^{\rm AFC}$ belongs to the admissible
  set of the PP optimization problem.
 \end{rmk}
  
\begin{rmk}
  In contrast to flux correction based on limiting, the OB-PP approach
  does not rely on the availability of a consistent low-order BP
  discretization for the given conservation law. For
  example, even the positivity-preserving
  exact solution of $\eqref{ibvp-pde}$ with
  $\mathbf{f}(u)=\mathbf{v}u$ may violate the upper bound
  $u^{\max}=1$ for a concentration field $u$ if the given
  velocity field $\mathbf{v}=\mathbf{v}(\mathbf{x})$ is not exactly
  divergence-free \cite{fctmax}. The PP algorithm can easily be configured
  to produce the best physics-compatible approximations using
  optimal flux control to preserve the global bounds and the
  conservation property. 
\end{rmk}

Similarly to limiter-based monolithic AFC schemes, optimal
flux control can also be applied at the level of the spatial
semi-discretization \eqref{semi-l}. Since the BP criterion
\eqref{pp-mono2} is satisfied for a family of closed-form limiters
\cite{timelim}, the feasible set of the  semi-discrete optimization
problem
\beq\label{pp-obj-mono}
\mbox{minimize}\ \frac{1}{2} (\dot u^* - \dot u^T)^\top M_C (\dot u^* - \dot u^T)
 + \frac{\mu}{2} (\dot u^*)^\top (M_L-M_C) \dot u^* 
 \qquad \mbox{subject to}
\eeq
\beq\label{pp-constr-mono}
c_i(u_i^{\min}- u_i)
\le r_i(u)+
\sum_{j\in\mathcal N_i\backslash\{i\}}m_{ij}(\dot u_i^*-\dot u_j^*)
\le c_i(u_i^{\max}- u_i)
\eeq
is nonempty, and an OB-PP algorithm can be used to fine-tune
the fluxes at individual RK stages. The potential advantages of flux
correction at the semi-discrete level include flexibility in the
choice of the time discretization and better convergence behavior
at steady state \cite{convex,timelim}.

\begin{rmk}
  For numerical schemes with diagonal mass matrices, PP
  optimization should be performed  using the coefficients
  $s_{ij}$ of another discrete Laplacian operator to define
  $g_{ij}^*=s_{ij}(\dot u_i^*-\dot u_j^*)$.
\end{rmk}

\section{Solution of optimization problems}
\label{sec:solver}

The PP optimization problem can be solved using a
barrier method, which guarantees that 
intermediate solutions stay in the feasible
set defined by \eqref{pp-constr}. Specifically,
we choose the primal-dual Newton method presented
in \cite[Section 6.6.2]{Ruszczynski+2011}.
The problem at hand can be written as
  $$ \min_{\dot{u} \in \Real^n} f(\dot u)\qquad
\text{s.t.} \quad A \dot u \le b, $$
where $n=N_h$ is the number of optimization variables. The
 inequality constraints are defined using
  $$ A = \begin{pmatrix} M_L - M_C \\ M_C - M_L \end{pmatrix} \in \Real^{2n \times n} \quad
  \text{and} \quad b = \begin{bmatrix} \frac{M_L}{\Delta t}(u^{\max} - u^n) - r \\ 
  \frac{M_L}{\Delta t}( u^n - u^{\min}) + r  \end{bmatrix} \in \Real^{2n}.$$
The gradient and Hessian of the convex objective function
$$ f_\mu(\dot u) = \frac{1}{2} (\dot u - \dot u^T)^\top M_C (\dot u - \dot u^T)
 + \frac{\mu}{2} \dot u^\top (M_L-M_C) \dot u $$
 are given by
   \begin{equation}\label{fmu_def}
  \nabla f_\mu = M_C (\dot u - \dot u^T) + \mu (M_L-M_C) \dot u
  \quad \text{and} \quad H_{f_\mu} = M_C + \mu (M_L-M_C).
  \end{equation}
  Introducing a vector  $s \in \Real^{2n}$ of slack variables, we
reformulate the inequality constraints as follows:
  $$ \min_{\dot{u} \in \Real^n} f_\mu(\dot u)\qquad
  \text{s.t.} \quad A \dot u + s = b,\qquad 
  s \ge 0.$$
The corresponding system of optimality conditions reads
\begin{align*}
 \nabla f_\mu( \dot u) + A^\top \lambda &= 0, \nonumber \\
 A \dot u + s &= b,  \label{opt_cond} \\
 s_i \lambda_i &= \sigma,\quad i=1,\ldots,2n, \nonumber
\end{align*}
where $\lambda\in\Real^{2n}$ is the vector of Lagrange multipliers and
$\sigma \ge 0$ is a parameter which is gradually decreased in the
process of solving a sequence of auxiliary problems.

  The initial guess $\dot u^{(0)}$ for the iterative optimization procedure
  should belong to the feasible set and be 
  a usable approximation to the optimization target $\dot u^T$.
  Such an approximation can be defined using an AFC scheme.
  Given the artificial diffusion operator $D$ and the
  antidiffusive fluxes $f_{ij}^*$ produced by
 a closed-form limiter, we
  compute $\dot u^{(0)}$ by solving the linear system
  \begin{equation}
    \label{udot_init}
    (M_L-M_C)\dot u^{(0)} = D u^n + f^*_{ij}
    \end{equation}
  subject to the equality constraint $\sum_{i=1}^{N_h}\dot u_i^{(0)}=0$
  which ensures uniqueness of the solution to the discrete Neumann
  problem with the singular graph Laplacian $M_L-M_C$. Unless
  mentioned otherwise, we calculate the fluxes $f^*_{ij}$
  for \eqref{udot_init} using
  Zalesak's FCT algorithm \cite{zalesak79}.
  
  We use Newton's method to update the solution $(\dot u,\lambda,s)$. At an intermediate step $j+1$, the
  search directions $\delta \dot u^{j+1}$ and $\delta\lambda^{j+1}$ are
  determined by solving a linear system of the form
 \begin{align*}
 \begin{pmatrix} H_f & A^T \\ A & G^j \end{pmatrix} \begin{bmatrix} \delta \dot u^{j+1} \\
 \delta\lambda^{j+1} \end{bmatrix} = \begin{bmatrix} - \nabla f_\mu(\dot u^j) - A^T \lambda^j \\
 b - A \dot u^j - \sigma / \lambda^j \end{bmatrix}.
 \end{align*}
 The block
 $G^j \in \Real^{2n \times 2n}$ is a diagonal matrix with entries $g^j_{ii} = - s^j_i / \lambda^j_i$; see \cite[Section 6.6.2]{Ruszczynski+2011}.
 
 The search directions $\delta s^{j+1}_i$ for individual
 slack variables $s_i$ are  defined by
 $$ \delta s_i^{j+1} = (\sigma - s^j_i \lambda^j_i - s^j_i \delta \lambda_i^{j+1}) / \lambda^j_i,\qquad i=1,\ldots,2n.$$
 Using the above search directions, the employed barrier method
 updates the solution as follows:
 \begin{align*}
   s^{j+1} &= s^j + \alpha \delta s^{j+1} > 0,\qquad
   \alpha \in \left[ 0 , 1 \right],\\
    \dot u^{j+1} &= \dot u^j + \alpha \delta \dot u^{j+1},\\
   \lambda^{j+1} &= \lambda^j + \beta \delta \lambda^{j+1} > 0,\qquad
   \beta \in \left[ 0 , 1 \right].
 \end{align*}

 To achieve fast convergence, it is essential to adjust the value
 of $\sigma$ adaptively. If the initial value of $\sigma$ is chosen
 too small, the method may converge slowly or fail to converge.
 Hence, a sufficiently large value of $\sigma$ should be used at
 the beginning. If 
 $\sigma$ turns out to be so large that the current Newton
 step would increase the value of $f_\mu$, this step should be
repeated with a smaller value of $\sigma$.
 
In our implementation, we iterate using a fixed value of
$\sigma$ until the improvement factor $$\frac{ f_\mu( \dot u^{j+1})
 }{f_\mu(\dot u^j)}$$ becomes close to or greater than $1$. Then
 we decrease $\sigma$ and restart the Newton iteration using
 the current solution as an initial guess. This process is repeated
 until a threshold value $\sigma^{\min}$ is reached. The final
 result is an approximate solution of the PP optimization 
problem \eqref{pp-obj},\eqref{pp-constr}.

We use the direct solver UMFPACK from the SuiteSparse library
\cite{SuSp} to solve linear systems. 
Since  $G^j$ depends on $s^j$ and $\lambda^j$,
the  LU decomposition needs to
be updated in every Newton iteration.

\section{Case studies and numerical examples}
\label{sec:num}

In this section, we apply the OB-PP algorithm to linear advection
problems and to an anisotropic diffusion equation. Stationary problems
of the form \eqref{bvp} are solved using time marching. For
comparison purposes, we present numerical results obtained with
AFC schemes based on Zalesak's FCT algorithm \cite{afc1,zalesak79}
and the monolithic convex limiting (MCL) strategy \cite{convex,entropyHO}.
All methods under investigation are implemented in the open-source
C++  finite element library MFEM \cite{mfem}. The target discretization
for PP optimal control and closed-form flux limiting is chosen
individually in each experiment. The objective function for all
optimization problems is defined using $\mu=0.01$.

\subsection{Solid body rotation}

We begin with a popular
solid body rotation test \cite{afc1,leveque}. The
 unsteady linear advection equation
$$
  \pd{u}{t}+\nabla\cdot(\mathbf{v}u)=0
$$
is solved in $\Omega=(0,1)^2$ using the solenoidal velocity field
${\bf v}(x,y)=(0.5-y,x-0.5)^T$.
The initial condition, as defined by LeVeque \cite{leveque},
is given by
$$
   u_0(x,y)=\begin{cases}
    u_0^{\rm hump}(x,y)
     &\text{if}\  \sqrt{(x - 0.25)^2 + (y - 0.5)^2}\le 0.15, \\
    u_0^{\rm cone}(x,y)
     &\text{if}\ \sqrt{(x - 0.5)^2 + (y - 0.25)^2}\le 0.15, \\
     1 &\text{if}\ \begin{cases}
       \left(\sqrt{(x - 0.5)^2 + (y - 0.75)^2}\le 0.15 \right) \wedge \\
       \left(|x - 0.5| \ge 0.025 \vee \ y\ge 0.85\right),
     \end{cases}\\
      0 &  \text{otherwise},
    \end{cases}
$$
where
\begin{align*}
u_0^{\rm hump}(x,y)&=    \frac14 + \frac14 \cos \left(
     \frac{\pi \sqrt{(x - 0.25)^2 + (y - 0.5)^2}}{0.15}\right),\\
u_0^{\rm cone}(x,y)&= 1-\frac{\sqrt{(x - 0.5)^2 + (y - 0.25)^2}}{0.15}.
\end{align*}
On the inflow boundary $\Gamma_D$ of $\Omega$, we prescribe the
homogeneous Dirichlet boundary condition
$$u(x,y)=0\qquad \mbox{for}\ (x,y)\in\Gamma_D.$$

As a target discretization, we use the 
fourth-order Taylor-Galerkin (TTG-4A)
method \cite{quartapelle,selmin_fef93}
  \begin{align*}
M_Cu^{n+1/3}&=M_Cu^n+\frac{\Delta t}{3}[Ku^n+b^n]
    +\frac{(\Delta t)^2}{12}Su^n,\\
M_Cu^T&=M_Cu^n+\Delta t[Ku^n+b^n]
  +\frac{(\Delta t)^2}{2}Su^{n+1/3},
  \end{align*}
  where $ K=(k_{ij})_{i,j=1}^{N_h}$ and $S=(s_{ij})_{i,j=1}^{N_h}$
  are sparse matrices with entries
  $$k_{ij}=-\sum_{e=1}^{E_h}
   \int_{K_e}  \varphi_i\nabla\cdot(\mathbf{v}\varphi_j)\dx,
  \qquad s_{ij}=-\sum_{e=1}^{E_h}
   \int_{K_e}(\mathbf{v}\cdot\nabla\varphi_i)
   (\mathbf{v}\cdot\nabla\varphi_j)\dx.$$
   The Dirichlet boundary conditions are taken into account using
    the vector $b=(b_{i})_{i}^{N_h}$ of surface integrals
$$
\qquad b_{i}=\sum_{e=1}^{E_h}\int_{K_e\cap\Gamma_D}
\varphi_i(u_D-u_h)|\mathbf{v}\cdot\mathbf{n}|\ds
=\tilde b_i+\sum_{j\in\mathcal N_i\backslash\{i\}}f_{ij}^b,
$$
$$\tilde b_{i}=\sum_{e=1}^{E_h}\int_{K_e\cap\Gamma_D}
\varphi_i(u_D-u_i)|\mathbf{v}\cdot\mathbf{n}|\ds,\qquad
f_{ij}^b=(u_i-u_j)\sum_{e=1}^{E_h}\int_{K_e\cap\Gamma_D}
\varphi_i\varphi_j|\mathbf{v}\cdot\mathbf{n}|\ds.$$
Flux correction of FCT and MCL type is performed at the second step using the
target fluxes
$$f_{ij}=m_{ij}(\dot u_i^T-\dot u_j^T)+d_{ij}(u_i^n-u_j^n)
  -\frac{\Delta t}{2}s_{ij}\Big(u_i^{n+1/3}-u_j^{n+1/3}\Big)+f_{ij}^b,$$
$$d_{ij}=\max\{-k_{ij},0,-k_{ji}\}$$
such that
$$
M_Lu^T=M_Lu^n+\Delta t[(K+D)u^n+\tilde b^{n}]+f,\qquad
f_i=\sum_{j\in\mathcal N_i\backslash\{i\}}f_{ij}.
$$
The potentials $\dot u_i^T$ for calculation of $f_{ij}$ and
definition of the OB-PP objective function are given by
$$
\dot u^T=M_C^{-1}\left[Ku^n+b^{n}+\frac{\Delta t}{2}Su^{n+1/3}\right].
$$
All methods under investigation lead to nonlinear flux-corrected
approximations of the form
$$
M_Lu^{n+1}=M_Lu^n+\Delta t\left[Ku^n+b^{n}+\frac{\Delta t}{2}Su^{n+1/3}
+g^*\right],\qquad
g_i^*=\sum_{j\in\mathcal N_i\backslash\{i\}}g_{ij}^*.
$$
The FCT and MCL schemes differ in the way to compute the
limited antidiffusive components $f_{ij}^*$ of
$$g_{ij}^*=f_{ij}^*+d_{ij}(u_j^n-u_i^n)
-\frac{\Delta t}{2}s_{ij}\Big(u_j^{n+1/3}-u_i^{n+1/3}\Big)-f_{ij}^b.$$
Zalesak's FCT limiter is designed to enforce the maximum principle
\eqref{lmp} for the solution of the fully discrete problem, while 
MCL constrains the spatial semi-discretization \eqref{semi-l} to
satisfy \eqref{pp-mono},\eqref{pp-mono2} with
$$
c_i=\sum_{j\in\mathcal N_i\backslash\{i\}}d_{ij}.
$$
The fully discrete OB-PP algorithm yields the optimal fluxes
$g_{ij}^*=m_{ij}(u_i^*-u_j^*)$ corresponding to the solution 
of \eqref{pp-obj},\eqref{pp-constr}. We do not consider the
semi-discrete version \eqref{pp-obj-mono},\eqref{pp-constr-mono}
in this example
because the selected target scheme corresponds to a specific
space-time discretization (TTG-4A).

\begin{figure}
  \subfloat[initial data, $u_h \in \begin{bmatrix} 0, 1 \end{bmatrix}$ \label{fig:3:ana}]{
\includegraphics[width=0.55\linewidth]{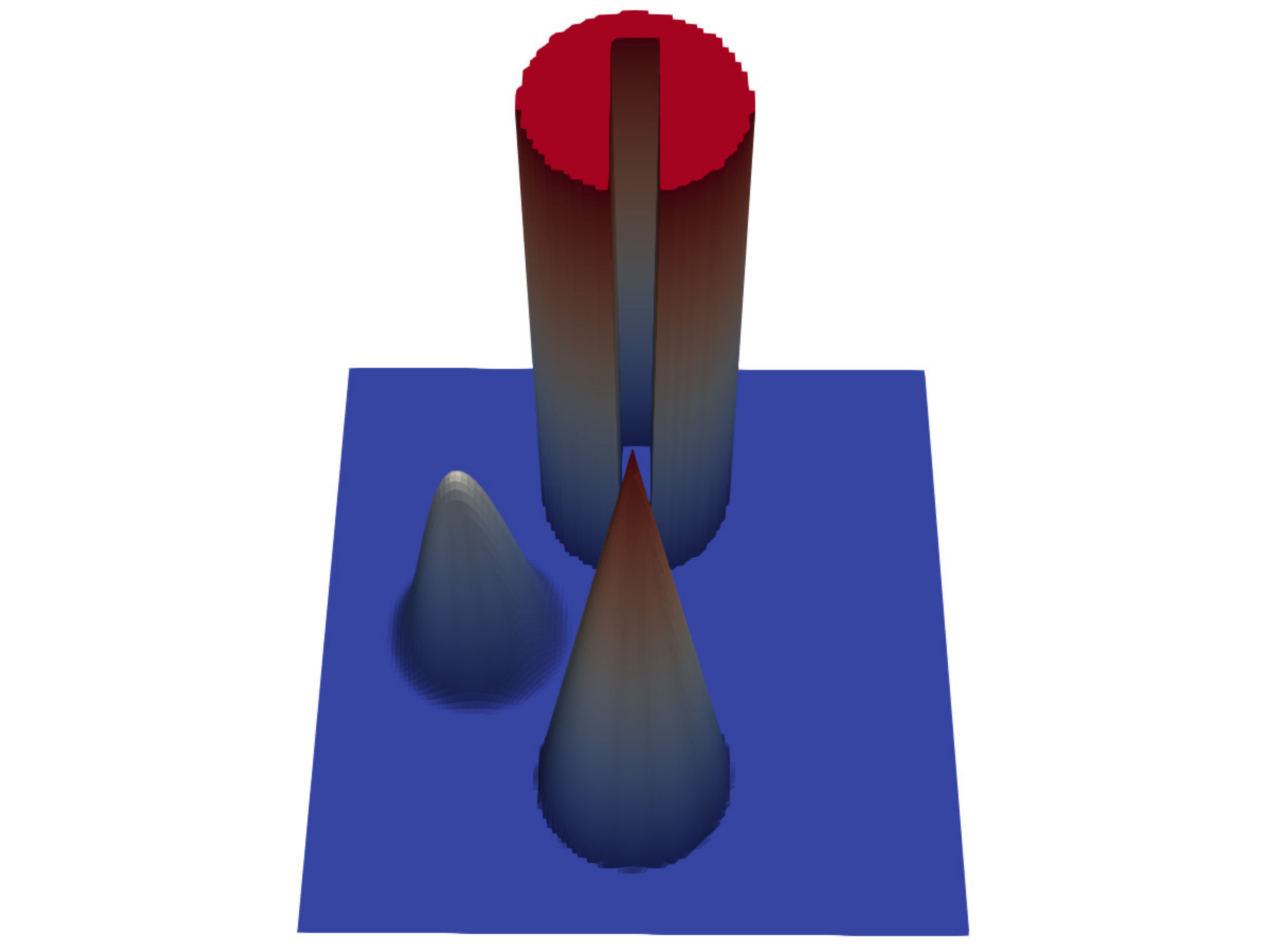}
}
\subfloat[MCL, $u_h \in \begin{bmatrix} 0, 0.98 \end{bmatrix}$ \label{fig:3:mcl}]{
\includegraphics[width=0.55\linewidth]{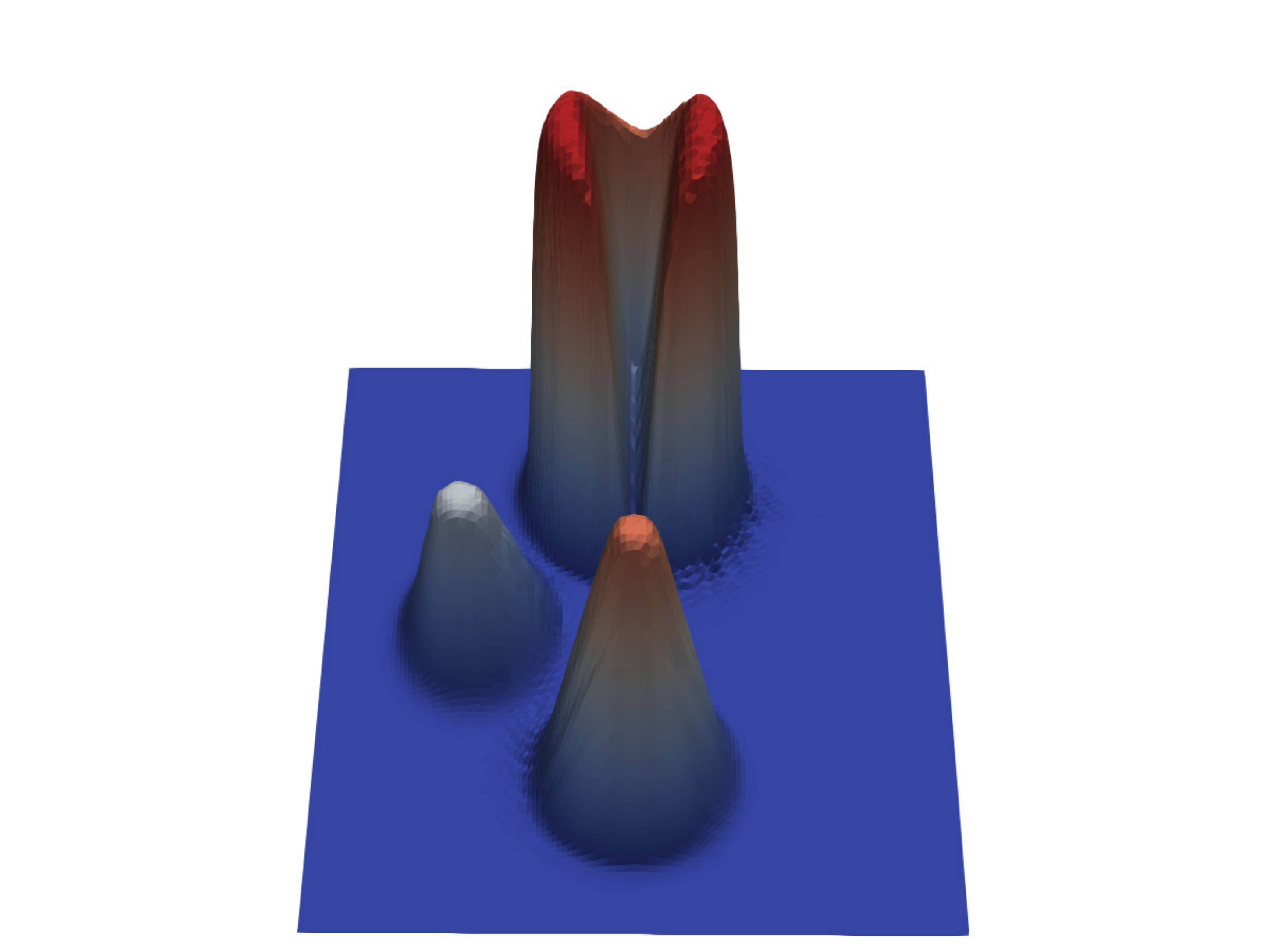}
}

\bigskip

\subfloat[FCT, $u_h \in \begin{bmatrix} 0, 0.988 \end{bmatrix}$ \label{fig:3:zal}]{
\includegraphics[width=0.55\linewidth]{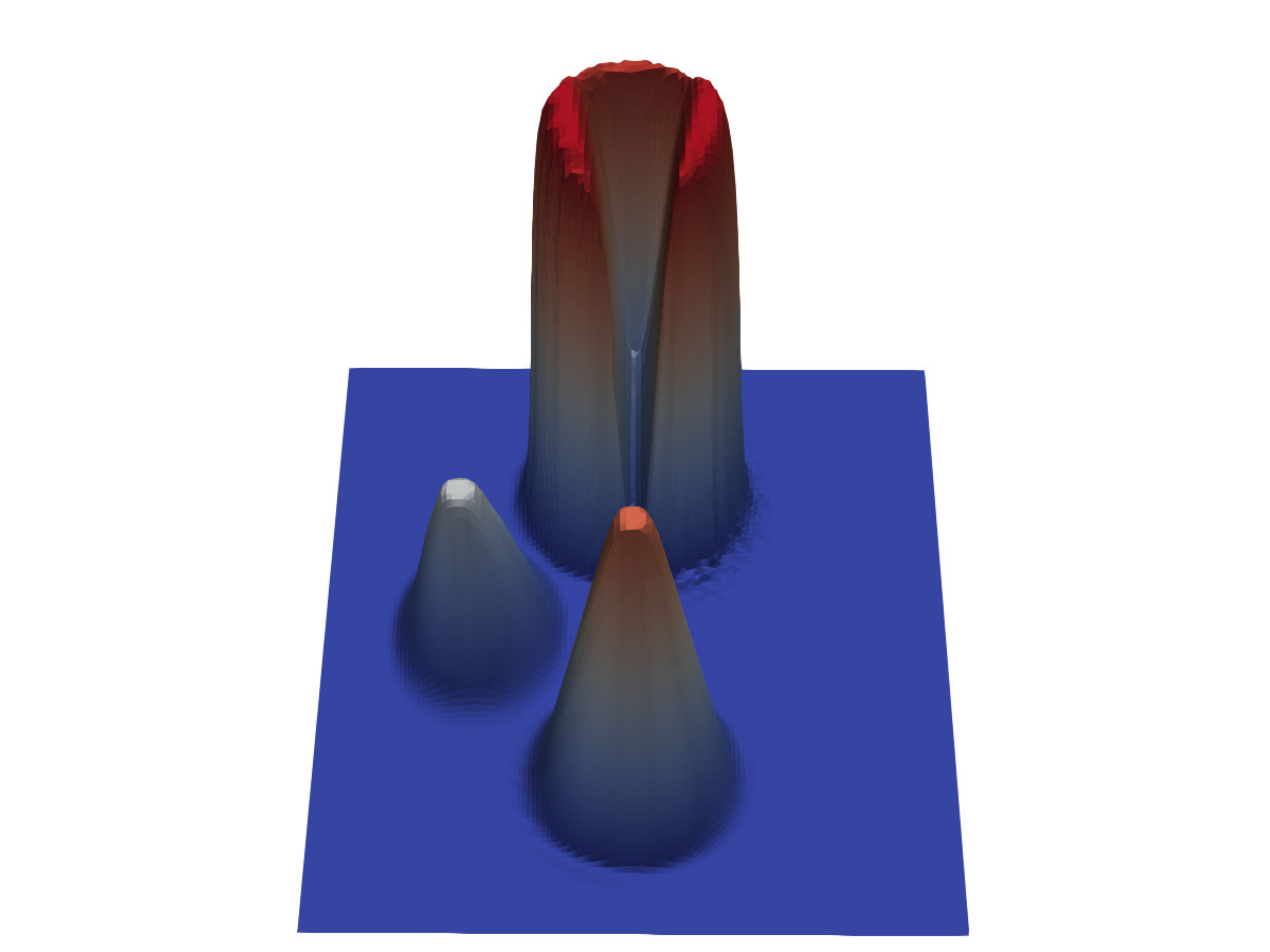}
}
\subfloat[OB-PP, $u_h \in \begin{bmatrix} 3.14 \cdot 10^{-9}, 1-10^{-7} \end{bmatrix}$ \label{fig:3:ob_pp}]{
\includegraphics[width=0.55\linewidth]{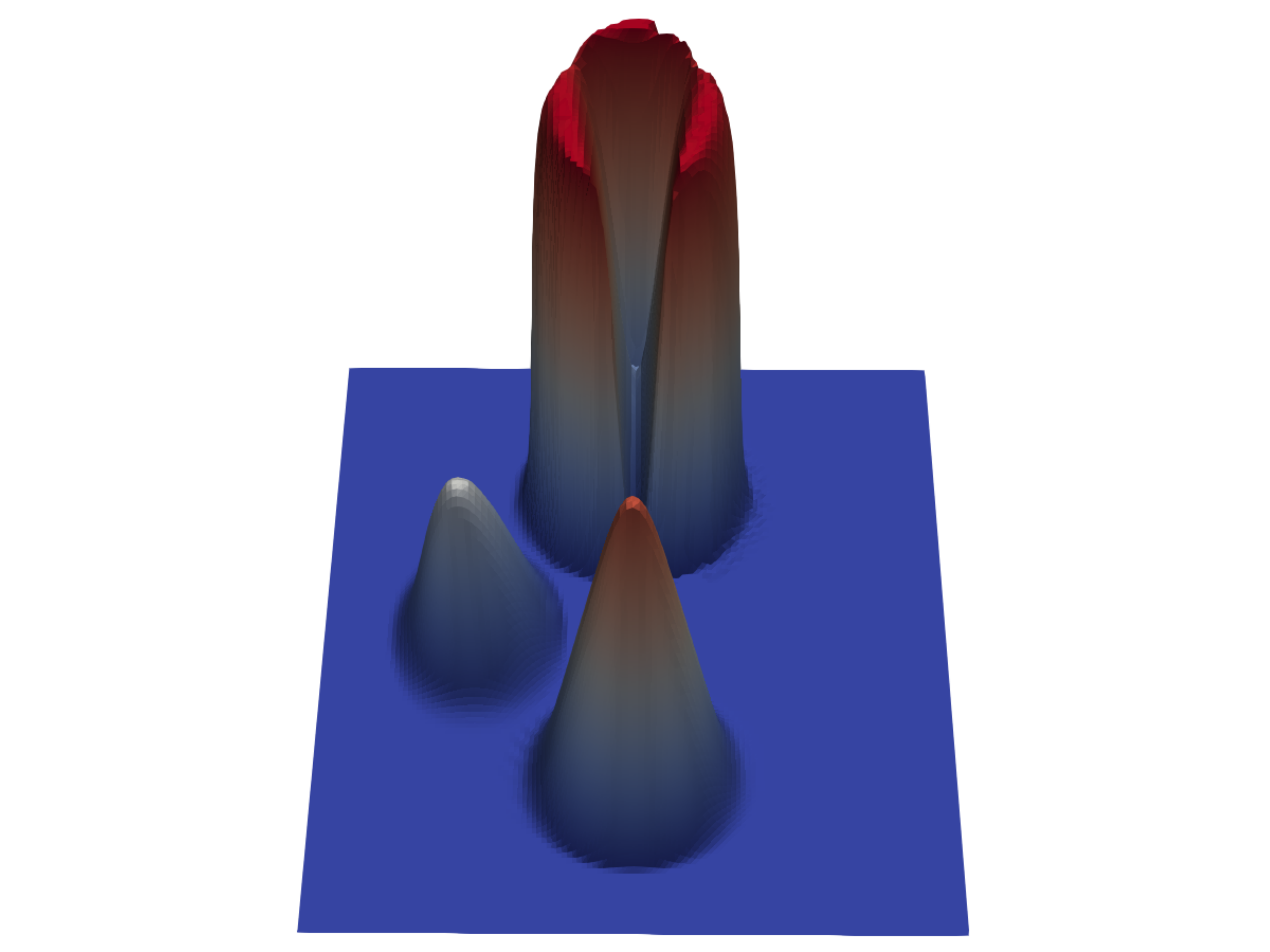}
}\vskip0.5cm

\caption{Solutions to the solid body rotation problem at $T=2\pi$.}
\label{fig:3:res_ss}
\end{figure}

In our numerical experiments for this benchmark,  we run simulations up
to the final time $T=2\pi$ using the mesh size $h = 1/128$ and
time step $\Delta t = 10^{-3}$. The interpolant $u_h(\cdot,0)$
of the initial data
is depicted in Fig.~\ref{fig:3:ana}. The
analytical solution $u_h(\cdot,2\pi)$
after one full rotation reproduces it exactly.
The numerical solutions shown in Figs~\ref{fig:3:mcl}
and  \ref{fig:3:zal} were obtained using closed-form flux
limiters of MCL and FCT type, respectively. The MCL result
exhibits significant levels of numerical diffusion on the narrow
back side of the slotted cylinder and at the two peaks. The
FCT solution is less diffusive but not as accurate as the
OB--PP result that we show in Fig.~\ref{fig:3:ob_pp}. Since
the optimization-based approach does not rely on worst-case
assumptions, the resulting approximations preserve the shape
of the initial data better their flux-limited counterparts.
In particular, the peak clipping effects are less pronounced,
and the global maximum of the OB-PP solution is closer to the
analytical value $u^{\max}=1$.

\subsection{Steady circular advection}

In the second standard test for AFC schemes, we solve the stationary
linear advection equation 
$$
\nabla\cdot(\mathbf{v}u)=0
$$
in $\Omega=(0,1)^2$ using the divergence-free velocity field
$\mathbf{v}(x,y)=(y,-x)$ unless mentioned otherwise. The inflow
boundary condition and the exact solution at any point
in $\bar\Omega$ are given by 
$$
u(x,y)=\left\{\begin{array}{ll}
1 &\quad \mbox{if} \ \ 0.15\le \sqrt{x^2+y^2}\le 0.45,\\
\cos^2\left(10\pi\frac{\sqrt{x^2+y^2}-0.7}{3}\right)
&\quad \mbox{if} \ \ 0.55\le \sqrt{x^2+y^2}\le 0.85,\\
0, & \quad\mbox{otherwise}.
\end{array}\right. $$

We march numerical solutions to the steady state using the
lumped-mass finite element version
$$
M_Lu^{n+1}=M_Lu^n+\Delta t[Ku^n+b(u^{n})]
  +\frac{(\Delta t)^2}{2}Su^{n}
  $$
 of the classical Lax-Wendroff (LW) method
 as target discretization for FCT, MCL, and OB-PP corrections. The
 amount of stabilizing  streamline diffusion is determined by the
 pseudo-time step $$\Delta t=\frac{\nu h}{\|\mathbf{v}\|_{L^\infty(\Omega)}},$$
 where $h$ is the mesh size and $\nu$ is a user-defined
 Courant number. The target  potential corresponding to the
 steady-state residual
 $r(u)=\left[K+\frac{\Delta t}{2}S\right]u+b(u)$ is given by $\dot u^T=0$.
 In our numerical study, we use $h=1/64$ and $\Delta t=10^{-3}$. For
 $\|\mathbf{v}\|_{L^\infty(\Omega)} = 1.413$, this choice of $\Delta t$
 corresponds to the Courant number $\nu = 0.09.$ 
  Computations are terminated at 
 $T=9.5$, when the $L^2$ norm of the pseudo-time derivative becomes smaller
 than $10^{-12}$ for MCL and as small as $7 \cdot 10^{-9}$ for OB-PP.
 The same output time is used for Zalesak's FCT algorithm which belongs
 to the family of predictor-corrector methods and, therefore, cannot
 be expected to produce a fully converged steady state solution.
 
\begin{figure}
\centering
\subfloat[exact, $u_h \in \begin{bmatrix} 0, 1 \end{bmatrix}$ \label{fig:1:ana}]{
\includegraphics[width=0.4\linewidth,trim=170 0 170 220,clip]{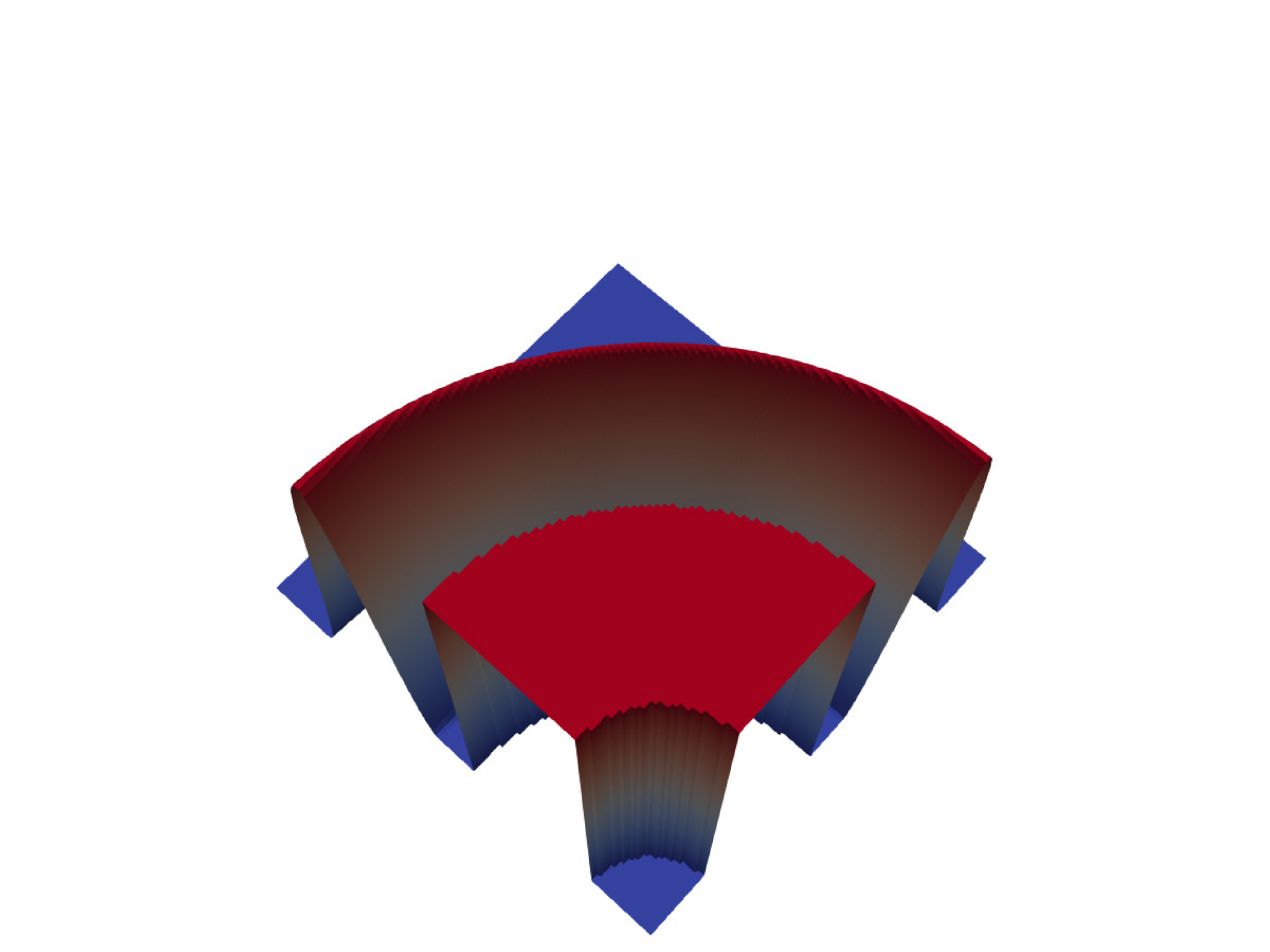}
}
\subfloat[unlimited, $u_h \in \begin{bmatrix} -0.24, 1.23 \end{bmatrix}$ \label{fig:1:no}]{
\includegraphics[width=0.4\linewidth,trim=170 0 170 220,clip]{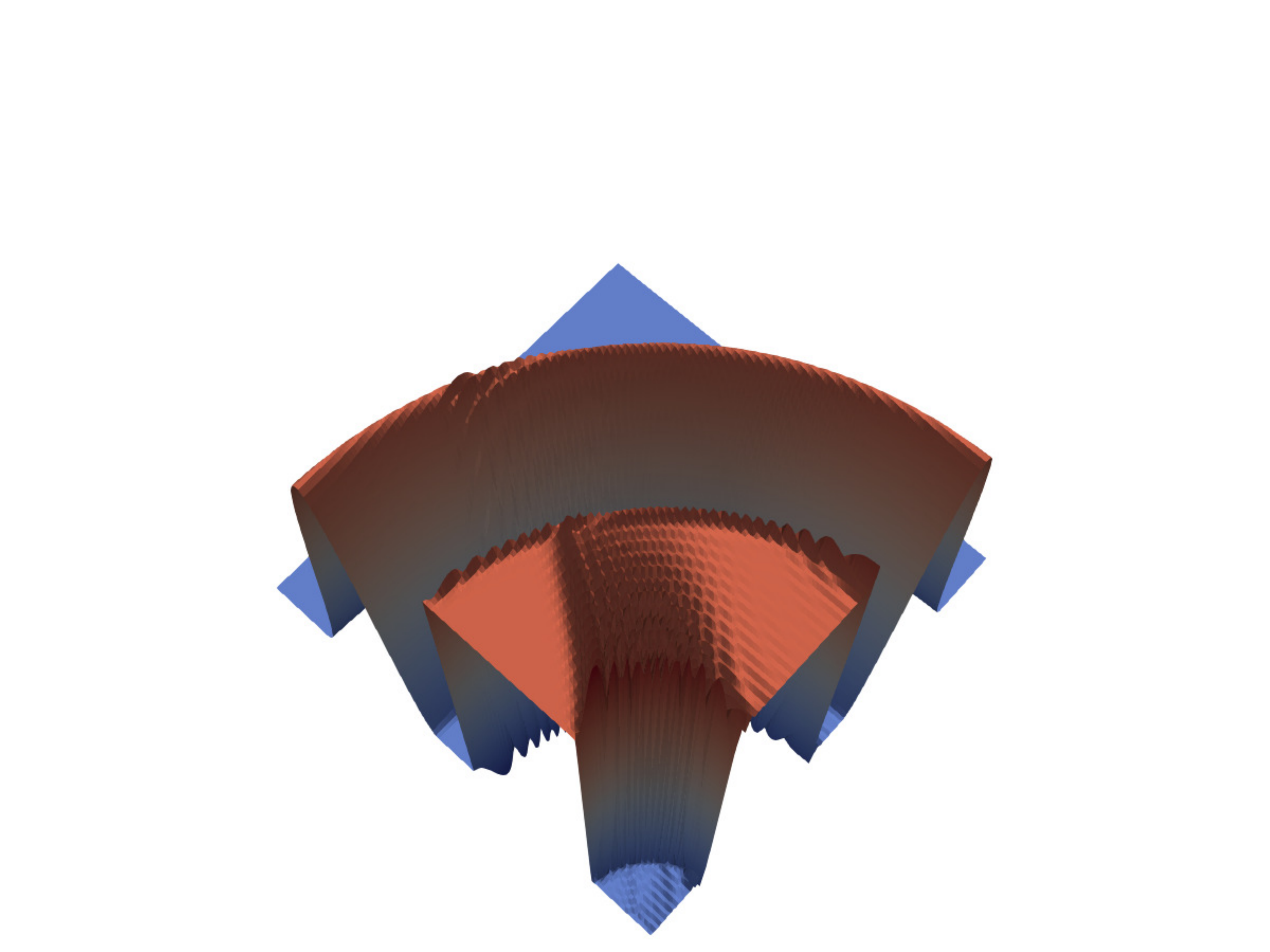}
}

\subfloat[MCL, $u_h \in \begin{bmatrix} 0, 1 \end{bmatrix}$ \label{fig:1:mcl}]{
\includegraphics[width=0.4\linewidth,trim=170 0 170 220,clip]{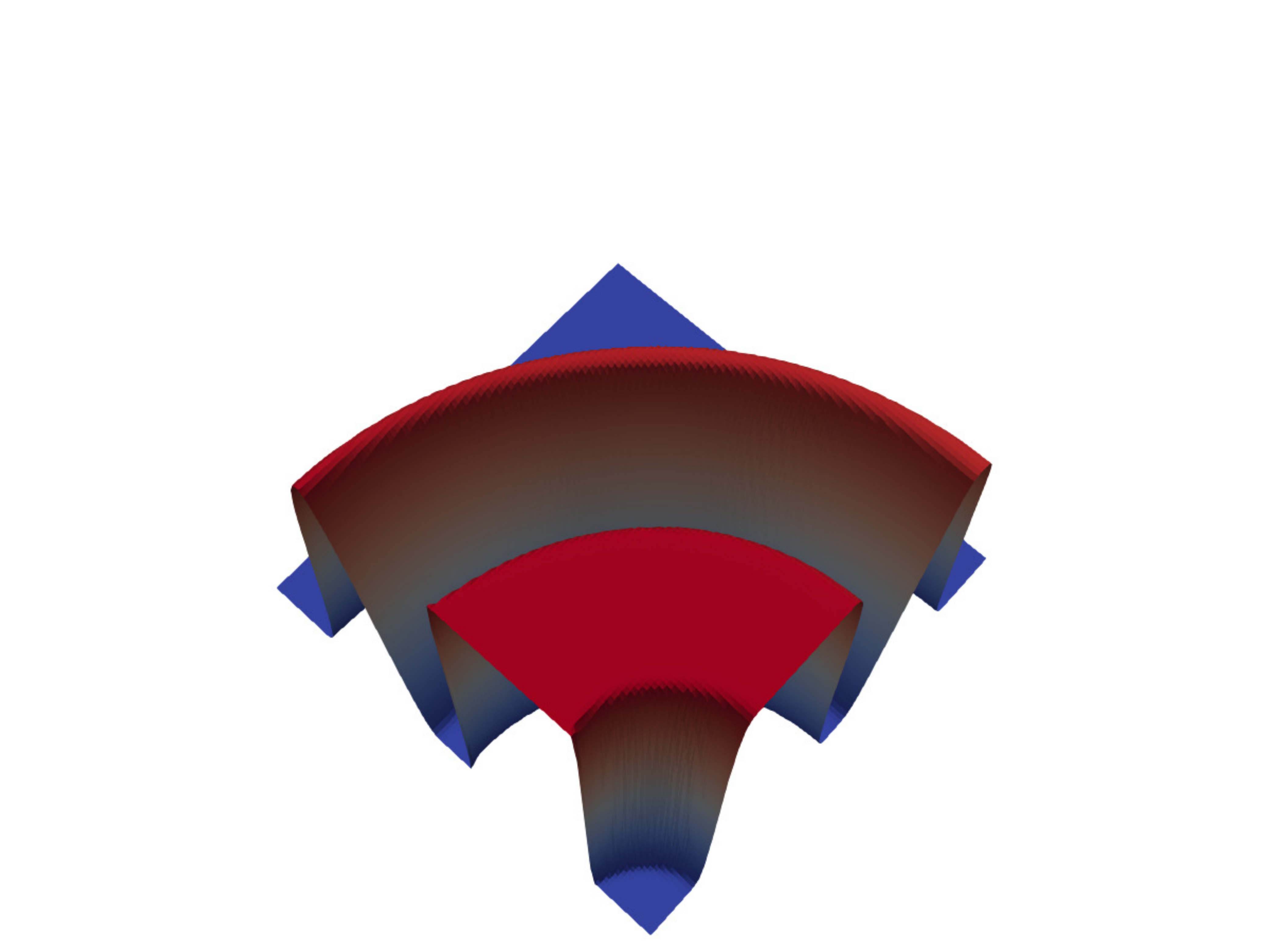}
}
\subfloat[FCT, $u_h \in \begin{bmatrix} 0, 1 \end{bmatrix}$ \label{fig:1:zal}]{
\includegraphics[width=0.4\linewidth,trim=170 0 170 220,clip]{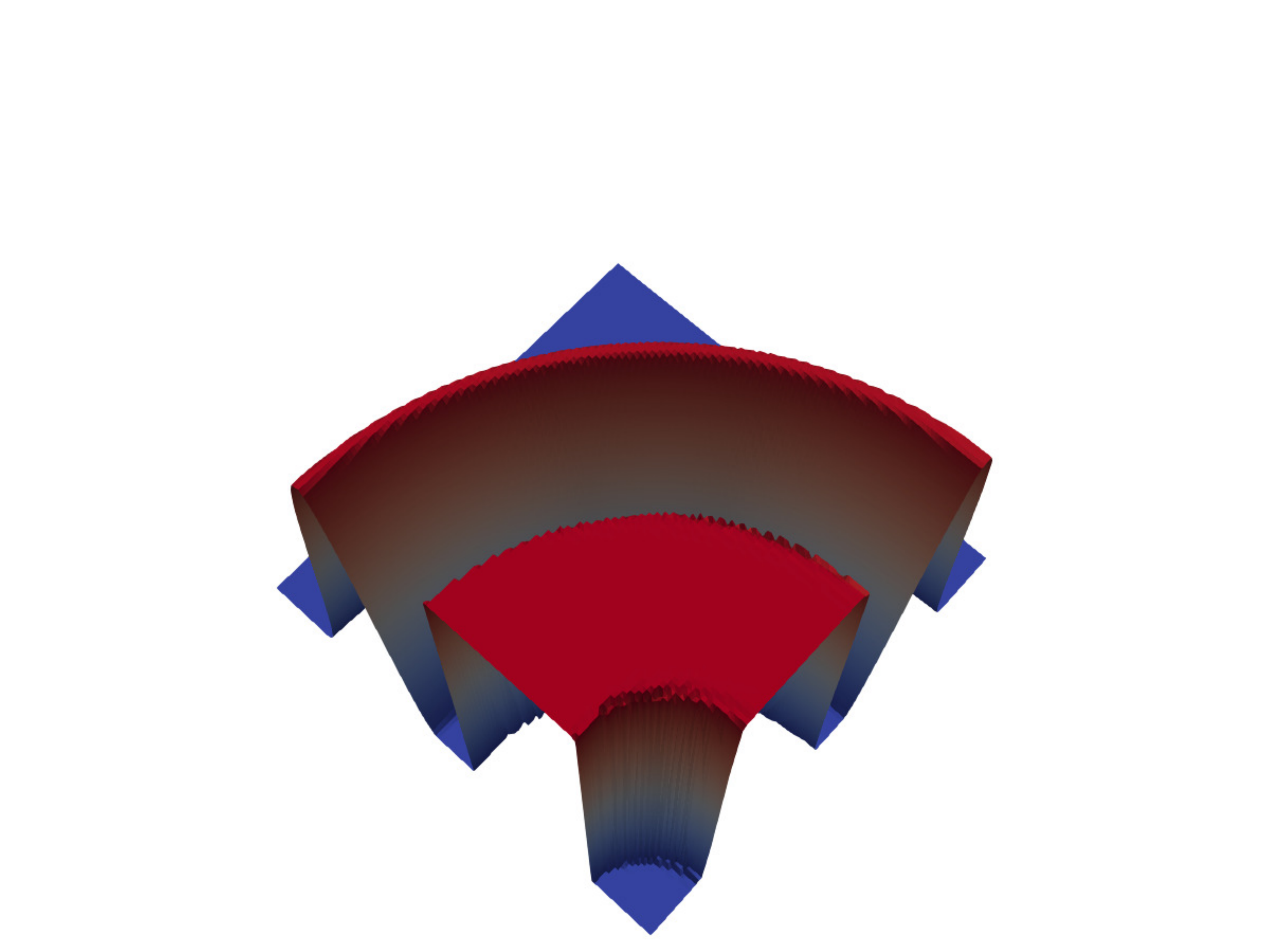}
}

\subfloat[OB-PP, $u_h \in \begin{bmatrix} 1.7 \cdot 10^{-11}, 1.0 \end{bmatrix}$ \label{fig:1:ob_pp}]{
\includegraphics[width=0.4\linewidth,trim=170 0 170 220,clip]{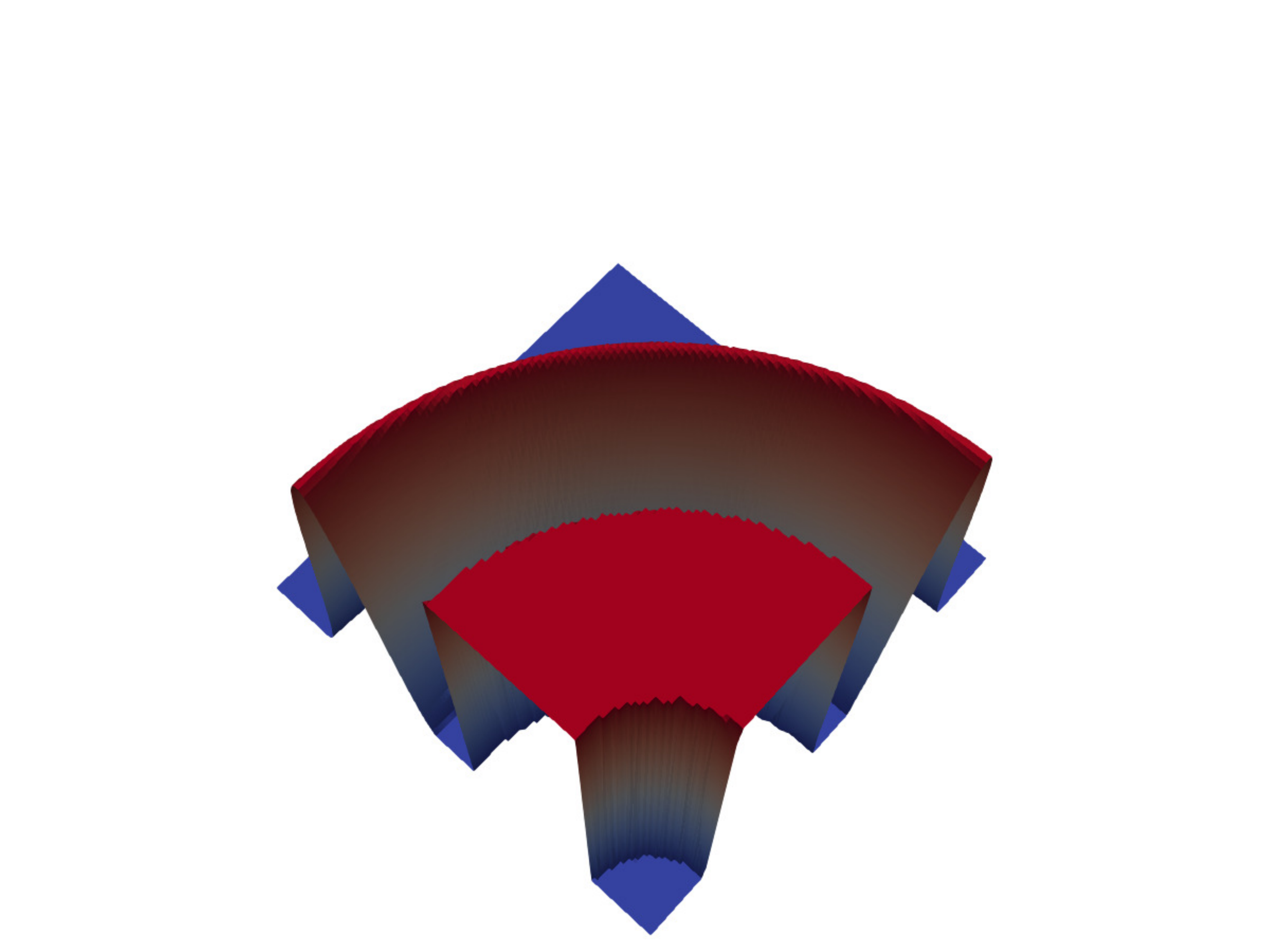}
}
\subfloat[$u_h^{\rm OB}-u_h^{\rm ex}\in \begin{bmatrix} -0.448, 0.975 \end{bmatrix}$ \label{fig:1:ob_diff}]{
\includegraphics[width=0.4\linewidth,trim=170 0 170 220,clip]{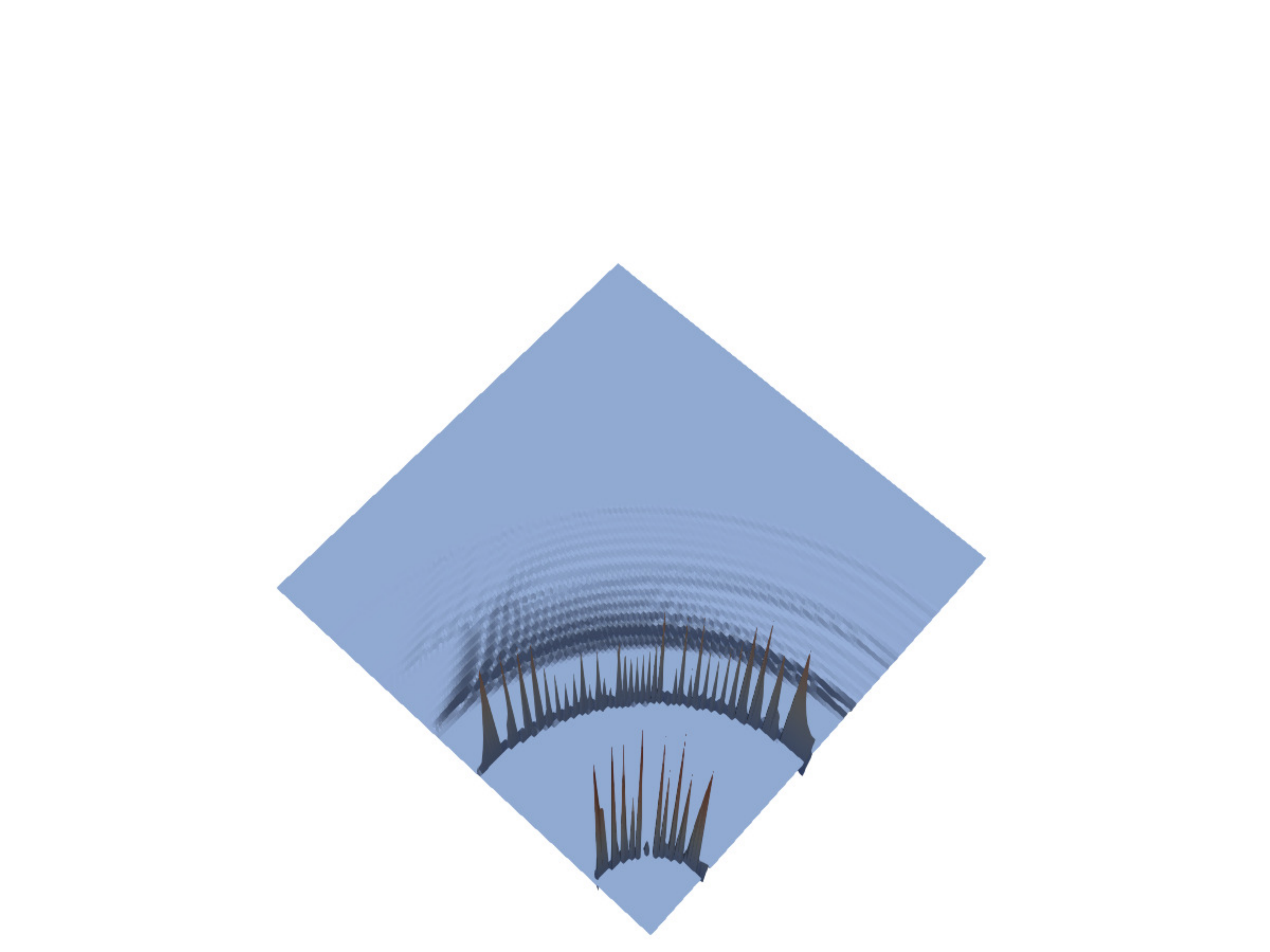}
}
\caption{Solutions to the steady advection problem at
 the pseudo-time $T=9.5$.}
\label{fig:1:res_ss}
\end{figure}

\begin{figure}
\centering
\includegraphics[width=0.4\linewidth,trim=50 0 50 50,clip]{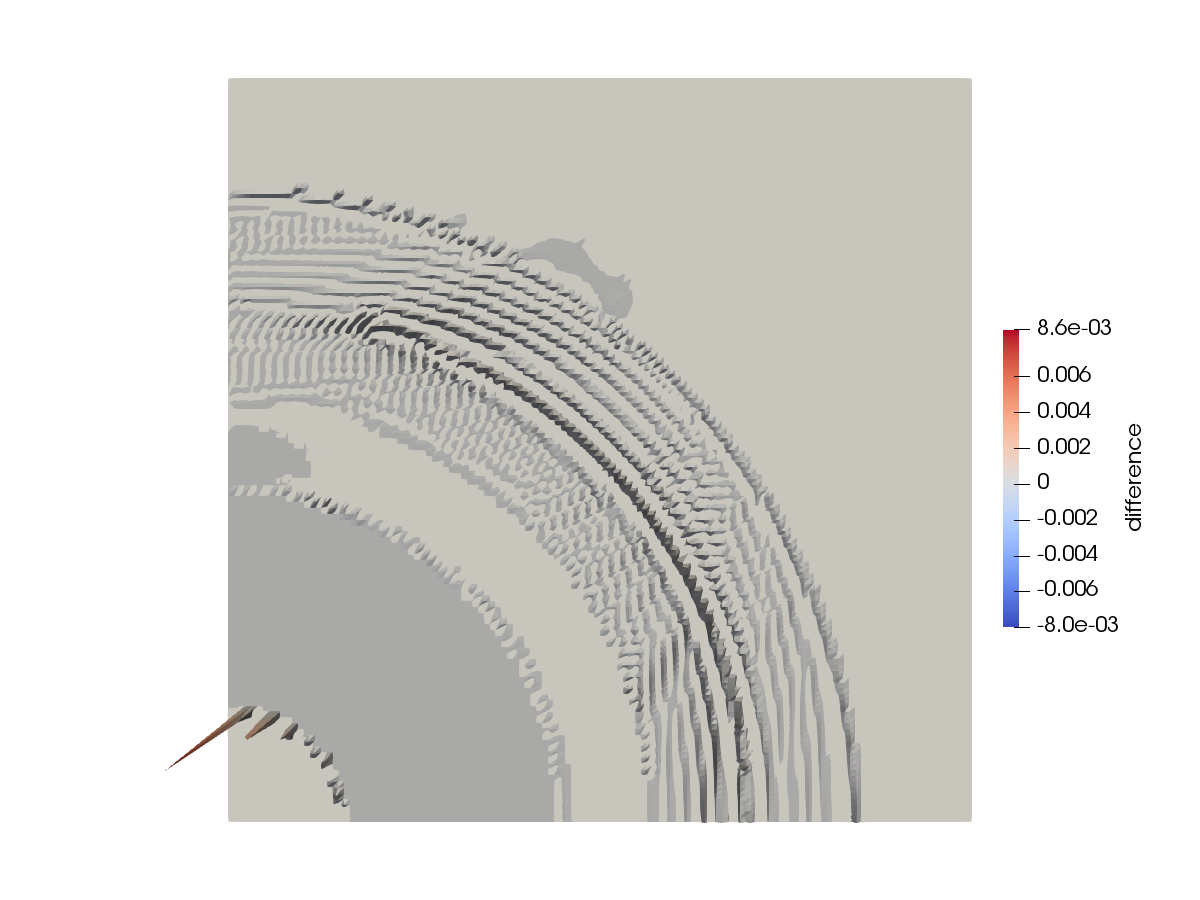}
\caption{Difference between the OB-PP results obtained
  using 
  formulations \eqref{pp-obj},\eqref{pp-constr} and
\eqref{pp-obj-mono},\eqref{pp-constr-mono}
  of the optimization problem
  (multiplied by a scaling factor of 50 for better
  visibility)\label{fig:1:diff}.}
\end{figure}

Figure \ref{fig:1:no} shows the unlimited Galerkin solution which
violates the discrete maximum principle and exhibits spurious oscillations
even in regions where the exact solution is smooth. All other
solutions shown in Fig.~\ref{fig:1:res_ss}
are globally bound preserving 
and  nonoscillatory. As in the unsteady advection
test, the most diffusive approximation is produced by MCL
(see Fig.~\ref{fig:1:mcl}), followed by FCT (cf. Fig.~\ref{fig:1:zal}). The
OB-PP result $u_h^{\rm OB}$ is virtually indistinguishable from the
interpolant $u_h^{\rm ex}$ of the
exact solution (compare Figs~\ref{fig:1:ana}
and \ref{fig:1:ob_pp}). To examine the distribution
of pointwise errors, we plot  $u_h^{\rm OB}-u_h^{\rm ex}$
in Fig.~\ref{fig:1:ob_diff}. As expected, the largest
errors are generated around the streamlines along which the
exact solution has a discontinuity, a peak, or a nonsmooth
transition in the crosswind direction. In this example,
we used OB-PP based on formulation \eqref{pp-obj},\eqref{pp-constr}.
The algorithm based on the semi-discrete version
\eqref{pp-obj-mono},\eqref{pp-constr-mono}
produces very similar results. The difference
between the solutions obtained with the two
OB-PP approaches is so small that we scale it by a
factor of 50 in Fig.~\ref{fig:1:diff} for better visibility.

In Table \ref{tab:i_zal_obpp}, we list the values of the
objective function \eqref{fmu_def} calculated using the
initial guess $\dot u^{(0)}$ (as defined by \eqref{udot_init}
with FCT fluxes $f_{ij}^*$)
and the final solution $\dot u^*$ of the OB-PP optimization
problem. It can be seen that the values of $f_\mu(\dot u^*)$ are
indeed significantly smaller than those of $f_\mu(\dot u^{(0)})$.

\begin{table}[h!]
\centering
\begin{tabular}{ccc}
\hline
pseudo-time & $f_\mu(\dot u^{(0)})$ &  $f_\mu(\dot u^*)$ \\ 
\hline
0.1 & 0.595 & 3.642e-02 \\
0.5 & 1.969 & 6.018e-02 \\
1.0 & 1.527 & 4.349e-02 \\
\hline
\end{tabular}
\caption{Values of the objective function $f_\mu$ for the
initial guess (Zalesak's FCT) and the OB-PP flux potential.}
\label{tab:i_zal_obpp}
\end{table}

\subsection{Anisotropic diffusion}

In the last example, we use the OB-PP algorithm to solve the steady anisotropic 
diffusion equation
$$
-\nabla\cdot({\cal D}\nabla u)=0
\label{diffeq}
$$
with \cite{afc1,anis}
$$
\mathcal D = R(-\theta) 
\left( \begin{array}{cc} 100 & 0 \\ 0 & 1 \end{array} \right)
{\cal R}(\theta),\qquad {\cal R}(\theta) =\left( \begin{array}{rr} \cos\theta & \sin\theta 
  \\ -\sin\theta & \cos\theta\\ \end{array} \right),
\qquad \theta = \frac{\pi}{6}$$
in the domain $\Omega  = (0, 1)^2 \backslash[4/9, 5/9]^2$. The 
outer and inner boundary of $\Omega$ are denoted
by $\Gamma_0$ and $\Gamma_1$, respectively.
The Dirichlet boundary condition for this test is given by
$$
 u(x,y) = \left\{\begin{array}{rl}
-1\  &\  \mbox{ if } (x,y)\in\Gamma_0,  \\
1\  &\ \mbox{ if } (x,y)\in\Gamma_1.
\end{array}\right.
 $$
 The standard Galerkin discretization yields the stiffness
 matrix $ K=(k_{ij})_{i,j=1}^{N_h}$ with entries
 $$
\qquad k_{ij}=-\sum_{e=1}^{E_h}
   \int_{K_e} \nabla\varphi_i\cdot(\mathcal D\nabla\varphi_j)\dx.
   $$
 Although the  unconstrained Galerkin solution is known to possess
 the best approximation property w.r.t. the energy norm, it may
 violate the global bounds $u^{\min}=-1$ and
 $u^{\max}=1$ if the mesh is nonuniform and/or the diffusion tensor
 $\mathcal D$ is highly anisotropic \cite{afc1,anis}.

 The design
 of BP flux limiters for elliptic problems is more difficult than for
 hyperbolic conservation laws because of the additional requirement
 that the limiting procedure be linearity preserving \cite{anis}. 
 The FCT and MCL schemes that we used in the first two examples
 are tailored for hyperbolic problems and do not ensure linearity
 preservation. Therefore, we constrain the Galerkin discretization
 of the anisotropic diffusion equation using the OB-PP algorithm.
 In this example, we use the target potential $\dot u^T=0$ and impose
 the Dirichlet boundary condition strongly (as an equality constraint).

\begin{figure}[h!]
 \centering
 \subfloat[\label{fig:4:diff}]{
 \includegraphics[width=0.5\linewidth]{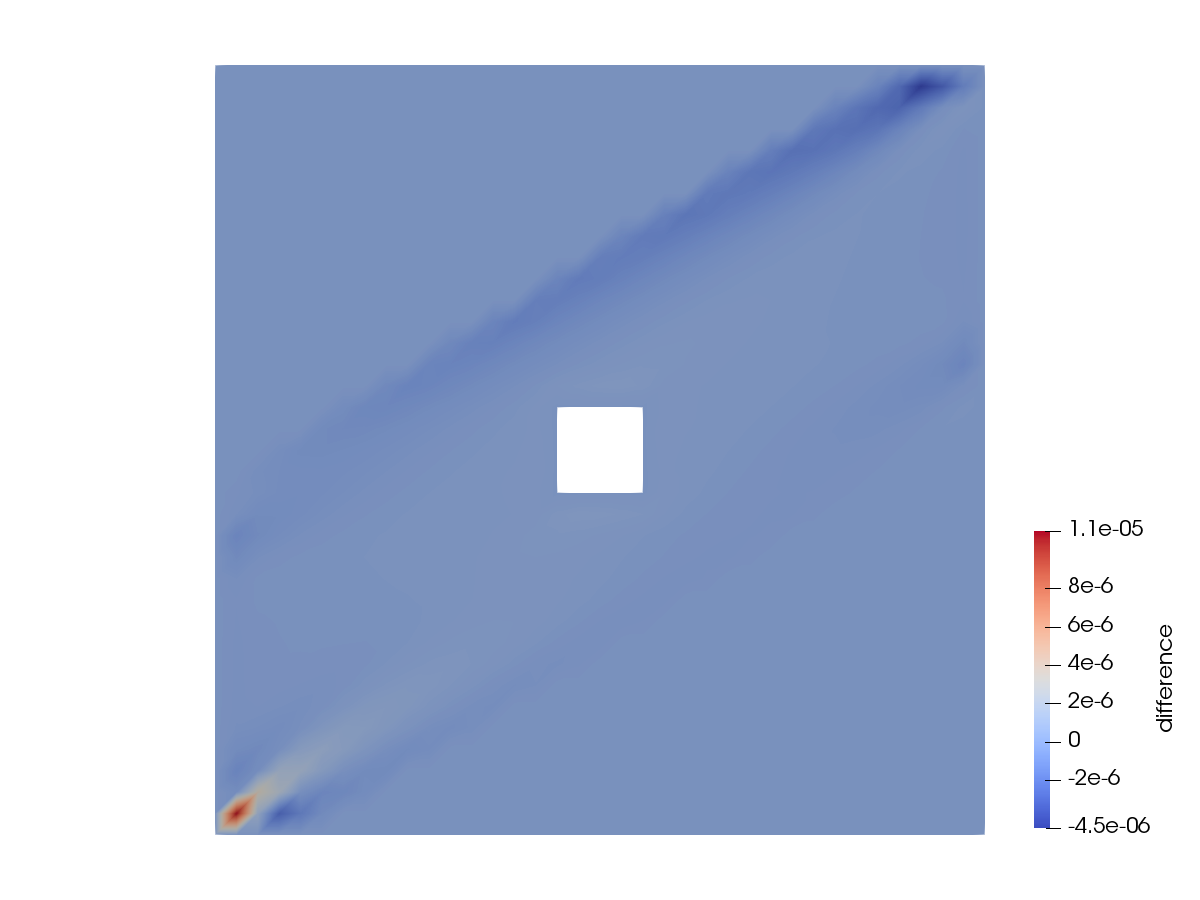}}
 \subfloat[]{
 \includegraphics[width=0.5\linewidth]{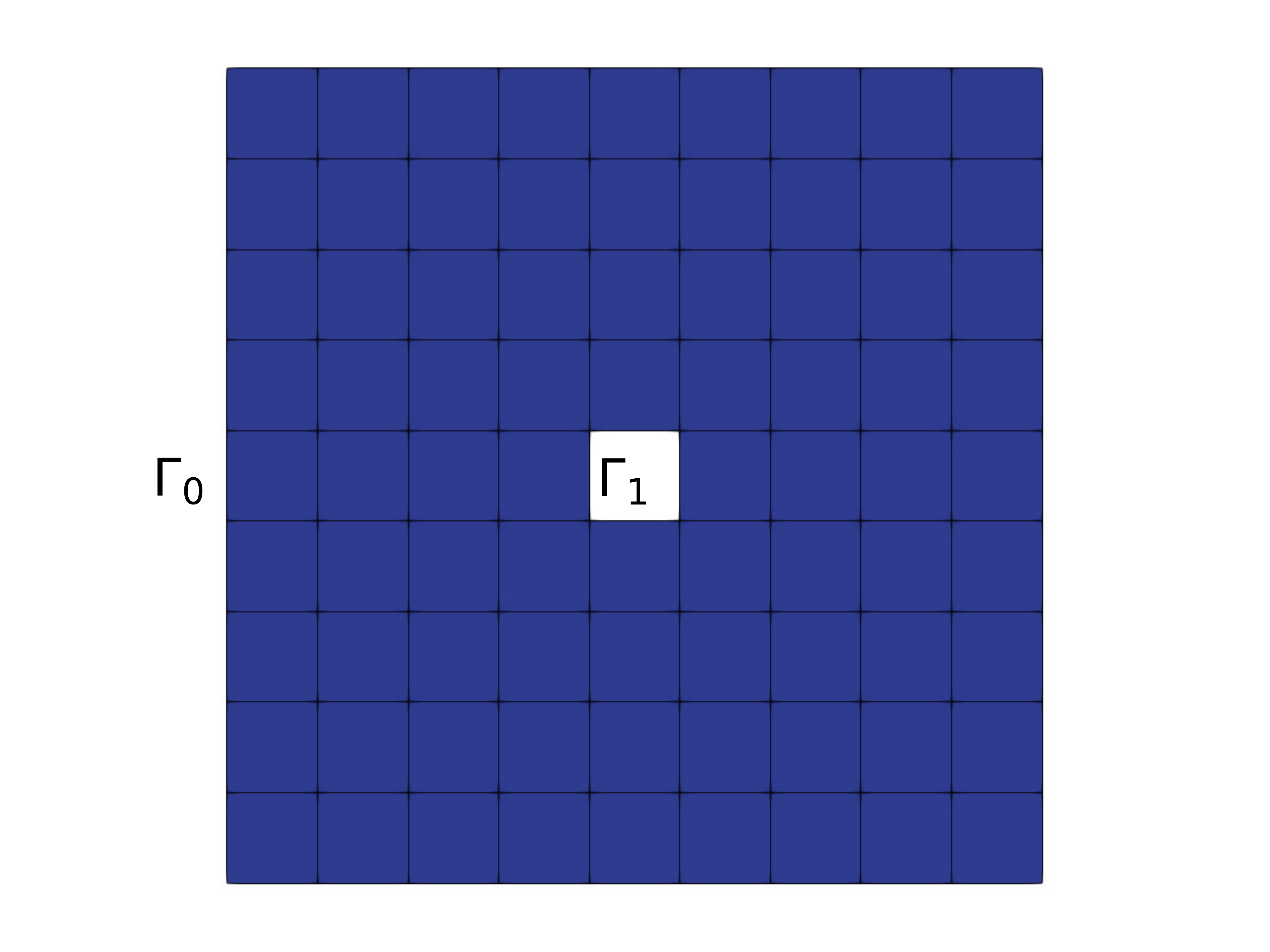}}
 \caption{Anisotropic diffusion test. (a) difference between the
 OB-PP solutions, (b) coarsest mesh.}
 \vskip0.5cm
 
 \subfloat[unlimited solution $u_h \in \begin{bmatrix} -1.004, 1 \end{bmatrix}$]{
 \includegraphics[width=0.48\linewidth]{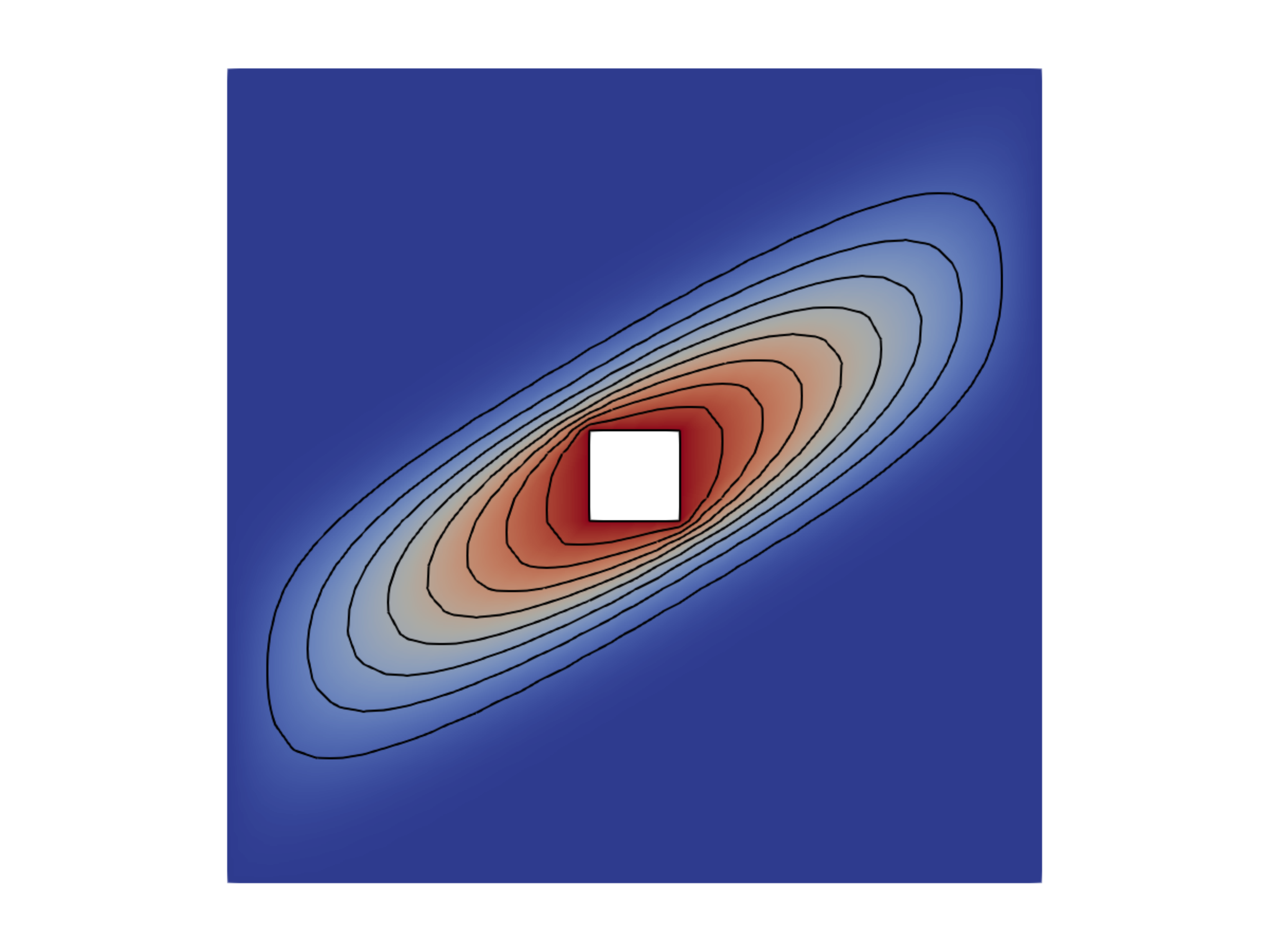}}
 \subfloat[OB-PP $u_h \in \begin{bmatrix} -1, 1 \end{bmatrix}$]{
 \includegraphics[width=0.48\linewidth]{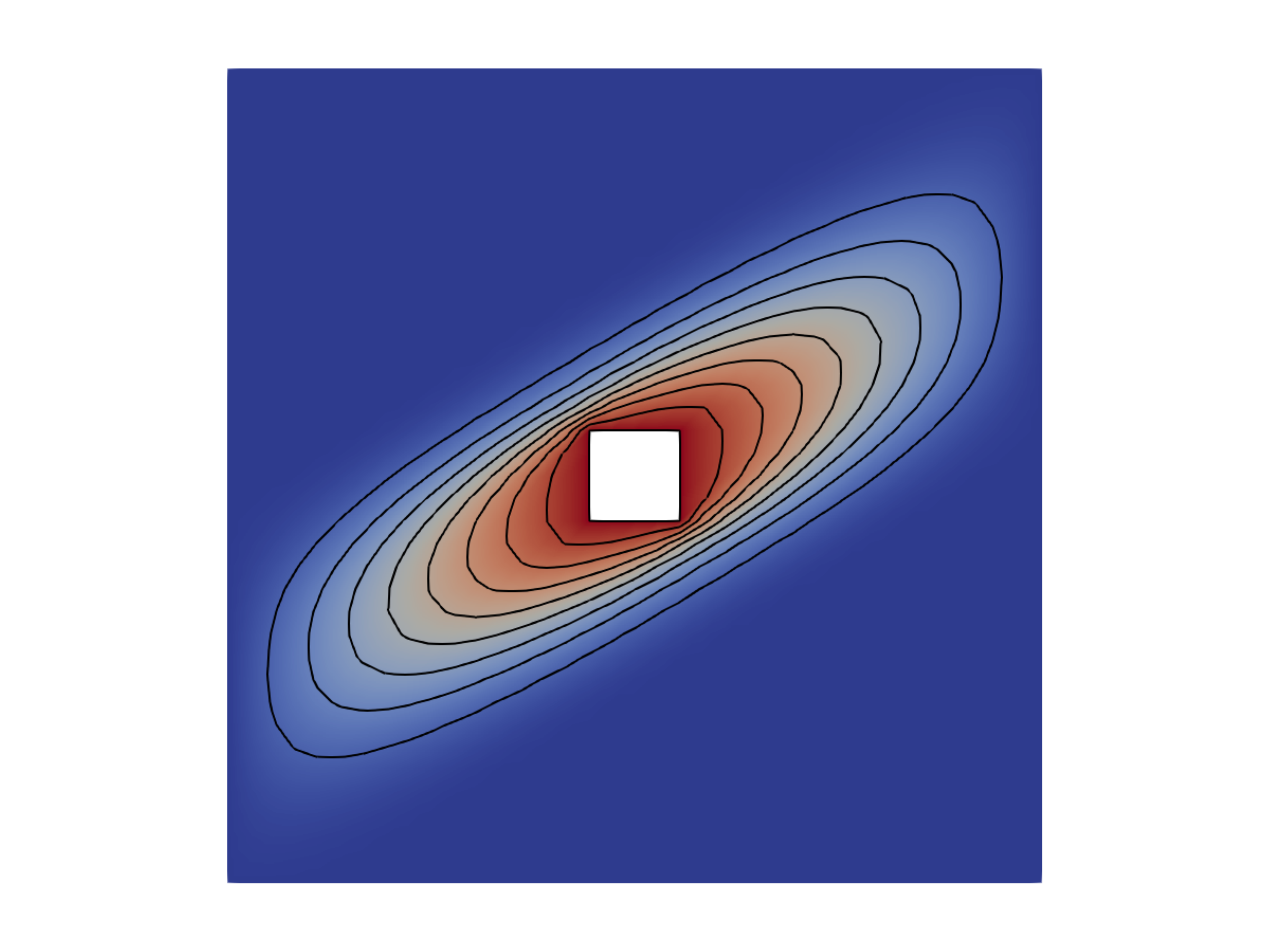}}
 \caption{Solutions to the anisotropic diffusion problem at the
   pseudo-time $T=2 \cdot 10^{-2}$.}
 \label{fig:4:res_ad}
 \end{figure}

 We run steady-state simulations using $h=1/18$ and $\Delta t = 10^{-6}$
 up to the pseudo-time $T=2 \cdot 10^{-2}$. In addition to the standard
 Galerkin method, we test the OB-PP algorithms based on
 \eqref{pp-obj},\eqref{pp-constr} and
 \eqref{pp-obj-mono},\eqref{pp-constr-mono}. The local bounds for
 the latter version of optimal control are defined using
 $$
 c_i=\sum_{j\in\mathcal N_i\backslash\{i\}}|k_{ij}|.
 $$
 The difference between the OB-PP solutions is shown in
 Fig.~\ref{fig:4:diff}. The maximum pointwise discrepancy is
 as small as $1.1 \cdot 10^{-5}$.
 Hence, the fully discrete and semi-discrete versions of OB-PP perform
 similarly again. In Fig.~\ref{fig:4:res_ad} we show the unlimited
 Galerkin solution and the result of fully discrete
 OB-PP optimization. The two diagrams look
 alike but the Galerkin approximation violates the 
 global lower
 bound $u^{\min}=-1$, while the OB-PP result is free of
 undershoots and overshoots.

 The results of grid convergence studies for the anisotropic
 diffusion problem are reported in Table~\ref{tab:AD_CR}.
 On the coarsest mesh, OB-PP is more accurate
 than the unlimited Galerkin method. On the next
 refinement level, the $L^1$ errors of the two
 approximations are almost the same, which explains
 the difference in the convergence rates of the two
 algorithms. On finer meshes, the $L^1$ error of the
 OB-PP discretization stays as small and goes to zero as fast
 as that of the underlying Galerkin scheme. However, the overhead
 cost of optimization-based flux correction is high since
 the discrete problem becomes nonlinear and many fixed-point
 iterations (pseudo-time steps) are needed to solve it.
 
 \begin{table}[h!]
 \centering
 \begin{tabular}{ccccc}
 \hline
 $\frac{1}{h}$ & unlimited & $p$ & OB-PP & $p$\\
 \hline 
 18  & 6.4831e-02 &        & 5.7464e-02 &        \\ 
 36  & 3.2154e-02 & 1.0117 & 3.2269e-02 & 0.8325 \\ 
 72  & 1.4789e-02 & 1.1205 & 1.4790e-02 & 1.1256 \\ 
 144 & 6.3533e-03 & 1.2189 & - & - \\
 \hline
 \end{tabular}
 \caption{$L^1$ errors and
   convergence rates for the anisotropic diffusion test.}
 \label{tab:AD_CR}
 \end{table}

\section{Conclusions}
\label{sec:conclusions}

The research presented in this paper indicates that optimal
control  based on the use of flux potentials provides
a versatile tool for enforcing discrete maximum principles
in a locally conservative manner. As demonstrated
in our numerical studies, the new OB-PP algorithm
yields more accurate approximations than closed-form
flux limiters and is better suited for preserving the range of
physically admissible states
in the presence of modeling errors.  Moreover, the number of flux
potentials is less than the number of fluxes that are used
as control variables in the OB-FF version.
While the cost of solving  inequality-constrained
optimization problems is considerable, it may be
comparable to or even lower than that of iterative flux
limiting in monolithic AFC schemes for stationary problems.
We envisage that entropy
stability conditions (cf. \cite{entropyCG,entropyHO}) and/or
additional problem-dependent constraints can be included in
the formulation of the PP optimization problem. Another
promising avenue for further research is  the design of
new monolithic flux control algorithms for spatial
semi-discretizations of the form \eqref{semi-l} using
criterion \eqref{pp-mono},\eqref{pp-mono2}
to formulate the inequality constraints.
It is hoped that further development and analysis of flux-based
optimal control approaches will make them an attractive alternative
to more traditional PDE-constrained optimization methods for
conservation laws.

\medskip
\paragraph{\bf Acknowledgments}

The second author would like to thank Dr. Pavel Bochev and Dr.
Denis Ridzal (Sandia National Laboratories) for very helpful
discussions of the potential-state potential-target flux
optimization method  (named so by Dr. Bochev)
at early stages of this work.

\end{document}